\documentclass[11pt,twoside]{article}
\usepackage{times}
\usepackage{amsmath,amssymb,amsthm}
\usepackage{enumerate}
\usepackage{cite}
\usepackage{mathrsfs,graphicx}
\allowdisplaybreaks

\newcommand{\no}{\nonumber}

\newcommand{\R}{\mathbb R}
\newcommand{\N}{\mathbb N}

\newcommand{\p}{\partial}
\newcommand{\ve}{\varepsilon}
\newcommand{\f}{\frac}
\newcommand{\la}{\lambda}

\newcommand{\g}{\gamma}

\newcommand{\ds}{\displaystyle}
\newcommand{\RN}[1]{\textup{\uppercase\expandafter{\romannumeral#1}}}

\def\defeq{\stackrel{\rm def}{=}}

\def\supp{\operatorname{supp}}
\def\ls{\lesssim}

\def\cE{{\mathcal E}}

\def\cQ{{\mathcal Q}}

\def\cV{{\mathcal V}}

\def\cX{{\mathcal X}}
\def\cY{{\mathcal Y}}

\def\grad{\nabla}

\def\virgp{\raise 2pt\hbox{,}}
\def\cdotpv{\raise 2pt\hbox{;}}

\def\C{\mathop{\mathbb C\kern 0pt}\nolimits}
\def\DD{\mathop{\mathbb D\kern 0pt}\nolimits}
\def\EE{\mathop{{\mathbb E \kern 0pt}}\nolimits}
\def\K{\mathop{\mathbb K\kern 0pt}\nolimits}
\def\N{\mathop{\mathbb N\kern 0pt}\nolimits}
\def\Q{\mathop{\mathbb Q\kern 0pt}\nolimits}
\def\R{\mathop{\mathbb R\kern 0pt}\nolimits}
\def\SS{\mathop{\mathbb S\kern 0pt}\nolimits}
\def\ZZ{\mathop{\mathbb Z\kern 0pt}\nolimits}
\def\TT{\mathop{\mathbb T\kern 0pt}\nolimits}
\def\P{\mathop{\mathbb P\kern 0pt}\nolimits}
\def\ds{\displaystyle}

\def\dive{\mathop{\rm div}}
\def\curl{\mathop{\rm curl}}

\def\supp{\mathop{\rm supp}\nolimits}

\def\no{\noindent}
\def\defeq{\stackrel{\rm def}{=}}
\def\p{\partial}
\def\ls{\lesssim}
\def\ve{\varepsilon}
\def\vv{\tilde v}
\def\QQ{\tilde Q}

\def\p{\partial}

\def\ve{\varepsilon}
\def\f{\frac}

\def\na{\nabla}
\def\la{\lambda}

\def\th{\theta}

\def\g{\gamma}

\def\inte#1{
\displaystyle\mathop{#1\kern0pt}^\circ }

\def\cE{{\mathcal E}}

\def\cQ{{\mathcal Q}}

\def\cV{{\mathcal V}}

\def\cX{{\mathcal X}}
\def\cY{{\mathcal Y}}

\def\grad{\nabla}

\def\virgp{\raise 2pt\hbox{,}}
\def\cdotpv{\raise 2pt\hbox{;}}

\def\C{\mathop{\mathbb C\kern 0pt}\nolimits}
\def\DD{\mathop{\mathbb D\kern 0pt}\nolimits}
\def\EE{\mathop{{\mathbb E \kern 0pt}}\nolimits}
\def\K{\mathop{\mathbb K\kern 0pt}\nolimits}
\def\N{\mathop{\mathbb N\kern 0pt}\nolimits}
\def\Q{\mathop{\mathbb Q\kern 0pt}\nolimits}
\def\R{\mathop{\mathbb R\kern 0pt}\nolimits}
\def\SS{\mathop{\mathbb S\kern 0pt}\nolimits}
\def\ZZ{\mathop{\mathbb Z\kern 0pt}\nolimits}
\def\TT{\mathop{\mathbb T\kern 0pt}\nolimits}
\def\P{\mathop{\mathbb P\kern 0pt}\nolimits}

\def\dive{\mathop{\rm div}}
\def\curl{\mathop{\rm curl}}

\def\supp{\mathop{\rm supp}\nolimits}

\def\no{\noindent}
\def\defeq{\stackrel{\rm def}{=}}
\def\p{\partial}
\def\ls{\lesssim}
\def\ve{\varepsilon}
\def\vv{\tilde v}
\def\QQ{\tilde Q}

\theoremstyle{plain}

\theoremstyle{definition}

\theoremstyle{remark}

\numberwithin{equation}{section}

\title{Global smooth axisymmetric solutions to 2D compressible Euler equations of Chaplygin gases with non-zero vorticity}

\author{Fei Hou \qquad Huicheng
  Yin$^{}$\footnote{Fei Hou (\texttt{houfeimath$@$gmail.com}) and
    Huicheng Yin (\texttt{huicheng$@$nju.edu.cn}, \texttt{05407$@$njnu.edu.cn}) were supported by
    the NSFC (No.~11571177, No. 11731007) and the Priority Academic Program
    Development of Jiangsu Higher Education Institutions.}\\
    [12pt] {\small School of Mathematical Sciences and Mathematical Institute, }\\
  {\small Nanjing
  Normal University, Nanjing 210023, China}}

\begin{document}

\date{}
\maketitle
\thispagestyle{empty}

\begin{abstract}
For 2D compressible isentropic Euler equations of polytropic gases, when the rotationally
invariant data are a perturbation of size $\ve>0$ of a rest state, S.~Alinhac in \cite{Alinhac92} and \cite{Alinhac93} establishes that the smooth solution blows up in finite time and  the lifespan $T_{\ve}$ satisfies $\ds\lim_{\ve\to 0}\ve^2 T_{\ve}=\tau_{0}^2>0$.
In the present paper, for 2D compressible isentropic Euler equations of Chaplygin gases,
we shall show that the small perturbed smooth solution exists globally when the rotationally
invariant data are a perturbation of size $\ve>0$ of a rest state.
Near the light cone, 2D Euler equations of Chaplygin gases can be transformed into a second order quasilinear wave
equation of potential, which satisfies both the first and the second null conditions. This will lead to that the corresponding
second order quasilinear wave equation admits a global smooth solution near the light cone (see \cite{Alinhac01}).
However,  away from the light cone, the hydrodynamical waves of 2D Chaplygin gases have no decay in time and strongly affect the related acoustical waves.
Thanks to introducing a nonlinear ODE and taking some delicate observations,
we can distinguish the fast decay part and non-decay part explicitly so that the global energy estimates with different weights can be derived by involved analysis.

\vskip 0.2 true cm
\noindent
\vskip 0.2 true cm
\noindent
\end{abstract}

\vskip 0.6 true cm
\tableofcontents
%\centerline{\bf Contents}
%\contentsline {section}{\numberline {1}Introduction}{2}
%\contentsline {section}{\numberline {2}Some Preliminaries}{7}
%\contentsline {section}{\numberline {3}On the analysis   of solutions near the light cone}{11}
%\contentsline {section}{\numberline {4}On the analysis   of solutions away from the light cone}{15}
%\contentsline {section}{\numberline {5}Elementary energy estimates of solutions}{25}
%\contentsline {section}{\numberline {6}Energy estimates  of solutions near the light cone}{27}
%\contentsline {section}{\numberline {7}Energy estimates  of solutions away from the light cone}{29}
%\contentsline {section}{\numberline {8}Proof of Theorem 1.1}{37}
%\contentsline {section}{\numberline {}References}{37}

\section{Introduction}
In this paper, we are concerned with the global existence of a smooth axisymmetric solution
to 2D compressible Euler equations of Chaplygin gases with non-zero vorticity. The initial data problem of 2D isentropic Euler equations is
\begin{equation}\label{Euler}
\left\{
\begin{aligned}
&\p_t\rho+div (\rho u)=0\qquad \qquad \qquad \qquad \qquad \qquad \qquad \qquad\text{(Conservation of mass)},\\
&\p_t(\rho u)+div (\rho u \otimes u) + \nabla P=0\qquad \qquad \qquad \qquad\qquad \text{(Conservation of momentum)},\\
&\rho(0,x)=\rho_0(x), u(0,x)=u_0(x),\\
\end{aligned}
\right.
\end{equation}
where $t\ge 0$, $x=(x_1, x_2)\in\Bbb R^2$, $\nabla=(\p_{x_1}, \p_{x_2})$,
and $u=(u_1,u_2),~\rho,~P$ stand for the velocity, density, pressure of gases respectively.
In addition, $\rho_0(x)>0$, $u_0(x)=(u_0^1(x), u_0^2(x))$,
$(\rho_0(x), u_0(x))\in C^{\infty}(\Bbb R^2)$, the pressure function $P=P(\rho)$ is smooth
and $P'(\rho)>0$ holds for $\rho>0$. When $P(\rho)=A\rho^{\g}$
with some positive constants $A$ and $\gamma$ ($1\le\gamma<3$), the corresponding compressible flows
are called the polytropic gases.

For the Chaplygin gases, the equation of pressure  state (one can see \cite{CF:book} and \cite{Godin07})
is given by
\begin{align}\label{ChaplyginGas}
P(\rho)=P_0-\ds\frac{B}{\rho},
\end{align}
where $P_0>0$ and $B>0$ are positive constants, and $P(\rho)>0$ for $\rho>0$. For convenience, we always assume
$P_0=2$ and $B=1$ in the whole paper.

If $(\rho, u)\in C^1$ is a solution of \eqref{Euler} with $\rho>0$, then \eqref{Euler}  admits
the following
equivalent form
\begin{equation}\label{EulerC1form}
\left\{
\begin{aligned}
&\p_t\rho+div(\rho u)=0,\\
&\p_tu+u\cdot\nabla u+\ds\frac{\nabla P}{\rho}=0,\\
&\rho(0,x)=\rho_0(x), u(0,x)=u_0(x).\\
\end{aligned}
\right.
\end{equation}

Let $\Omega=: x_1\p_{2}-x_2\p_{1}$ be the rotation operator. Set $\tilde\Omega u_0(x)=:(\Omega u^1_0(x),\Omega u^2_0(x))+(u^2_0(x),-u^1_0(x))$.
Through the paper, we assume $\Omega\rho_0=0$ and $\tilde\Omega u_0=0$.
In this case, one easily knows that
$\rho_0(x)=\rho_0(r)$ and $u_0(x)=f_0(r)\frac{x}{r}+g_0(r)\frac{x^\perp}{r}$ with $x^\perp=(-x_2,x_1)$
and $r=\sqrt{x_1^2+x_2^2}$.
Then it follows from \eqref{EulerC1form} that the smooth solution $(\rho, u)$ will admit such a form for any
$t\ge 0$: $\rho(t,x)=\rho(t,r)$ and $u(t,x)=f(t,r)\frac{x}{r}+g(t,r)\frac{x^\perp}{r}$.
Denote by $v=:\frac{1}{\rho}-1$ and  $v_0(r)=:\frac{1}{\rho_0(r)}-1$. Then \eqref{EulerC1form}
together with \eqref{ChaplyginGas} is equivalent to the following problem for $(v, f, g)$:
\begin{equation}\label{axisymEuler}
\left\{
\begin{aligned}
    &\p_tv-(\p_rf+\frac1rf)=Q_1,\\
    &\p_tf-\p_rv=Q_2,\\
    &\p_tg+f\p_rg+\frac1rfg=0,\\
    &v(0,r)=v_0(r),~f(0,r)=f_0(r),
    ~g(0,r)=g_0(r),
\end{aligned}
\right.
\end{equation}
where the quadratic nonlinear terms $Q_1$ and $Q_2$ are
\begin{eqnarray*}
    &&Q_1=: v(\p_rf+\frac1rf)-f\p_rv,\label{Q1}\\
    &&Q_2=: v\p_rv-f\p_rf+\frac1rg^2.\label{Q2}
\end{eqnarray*}

For fixed integer $N\ge N_0=7$, denote the size $\ve$ of the perturbed initial data of \eqref{EulerC1form} around the rest state
$(1, 0, 0)$ by
\begin{align*}
&\ve=:\sum_{k=0}^N\|(r\p_r)^k(\rho_0-1,u_0)\|_{L^2}
+\sum_{k+l=0}^{N-2}\|(r\p_r)^k\grad^l\curl u_0\|_{L^3}\no\\
&\qquad+\sum_{k+l=0}^{N-1}\|<r>(r\p_r)^k\na^l
(\nabla\rho_0,\dive u_0,\curl u_0)\|_{L^2},
\end{align*}
where $<r>=(1+r^2)^\frac12\ge1$, and $\|h\|_{L^p}=:(\int_{\R^2}|h|^pdx)^{\f1p}$ for $p\ge 1$.
The main result in this paper is:

\vskip 0.2 true cm

{\bf Theorem 1.1.} {\it There exists a constant $\ve_0>0$ such that if the initial
data $(\rho_0(x),u_0(x))\in C^{\infty}(\Bbb R^2)$ satisfy $\Omega\rho_0=0$, $\tilde\Omega u_0=0$ and
$\supp(\rho_0-1,u_0)\subset\{r\le\frac18\}$, then for $\ve\le\ve_0$, problem \eqref{EulerC1form} together with \eqref{ChaplyginGas}
admits a global smooth solution $(\rho(t,x),u(t,x))$  which
fulfills $\Omega\rho(t,x)=0$ and $\tilde\Omega u(t,x)=0$.}

\vskip 0.2 true cm

{\bf Remark 1.1} {\it It is only for notational convenience that we assume $\supp(\rho_0-1,u_0)\subset\{r\le\frac18\}$ in Theorem 1.1.
In fact, for $\supp(\rho_0-1,u_0)\subset\{r\le M_0\}$ with any fixed number $M_0>0$, Theorem 1.1 still holds.}

\vskip 0.2 true cm

{\bf Remark 1.2} {\it For the polytropic gases with $P(\rho)=\frac1\gamma\rho^\gamma$ ($\gamma\ge 1$)
(without loss of generality and for simplicity, $A=\f{1}{\g}$ is assumed in the state equation $P(\rho)=A\rho^{\g}$) and for the symmetric solution
$(\rho(t,r), f(t,r)\frac{x}{r}+g(t,r)\frac{x^\perp}{r})$ of \eqref{EulerC1form},
by setting $c(\rho)=\sqrt{P'(\rho)}=1+\frac{\gamma-1}{2}\dot c(\rho)$ for $\gamma\neq1$
and $\dot c(\rho)=\ln\rho$ for $\gamma=1$, then it follows from the first and the second equation of \eqref{EulerC1form} that
\begin{equation}\label{symmetrized}
\left\{
\begin{aligned}
    &(\p_t+f\p_r)\dot c(\rho)+[1+\frac{\gamma-1}{2}\dot c(\rho)](\p_rf+\frac1rf)=0,\\
    &(\p_t+f\p_r)f+[1+\frac{\gamma-1}{2}\dot c(\rho)]\p_r\dot c(\rho)=\frac{g^2}{r}.
\end{aligned}
\right.
\end{equation}
Denote by $Z_\pm=: f\pm\dot c(\rho)$. Then \eqref{symmetrized} can be written as
\begin{equation}\label{diagonalized}
\left\{
\begin{aligned}
    &(\p_t+\p_r)Z_++\frac{\gamma+1}{4}Z_+\p_rZ_+=\frac{\gamma-3}{4}Z_-\p_rZ_+-\frac1rf[1+\frac{\gamma-1}{2}\dot c(\rho)]+\frac{g^2}{r},\\
    &(\p_t-\p_r)Z_-+\frac{\gamma+1}{4}Z_-\p_rZ_-=\frac{\gamma-3}{4}Z_+\p_rZ_-+\frac1rf[1+\frac{\gamma-1}{2}\dot c(\rho)]+\frac{g^2}{r}.
\end{aligned}
\right.
\end{equation}
From this, one knows that the main parts in \eqref{diagonalized} are both Burgers equations of $Z_+$ and $Z_-$.
Therefore, S.~Alinhac in \cite{Alinhac93}
establishes that the smooth solution of \eqref{EulerC1form} together with $P(\rho)=\frac1\gamma\rho^{\gamma}$, $\Omega\rho_0=0$
and $\tilde\Omega u_0=0$ ($\gamma\ge 1$) blows up in finite time.
}
\vskip 0.2 true cm
{\bf Remark 1.3} {\it For the Chaplygin gases with $P(\rho)=2-\frac1{\rho}$ and for the symmetric solution
$(\rho(t,r), f(t,r)\frac{x}{r}$ $+g(t,r)\frac{x^\perp}{r})$ of \eqref{EulerC1form}, one can rewrite the nonlinear terms $Q_1$ and $Q_2$
in \eqref{axisymEuler}  as
\begin{eqnarray*}
    &Q_1=v\p_r(v+f)-(v+f)\p_rv+\frac1rvf,\\
    &Q_2=v\p_r(v+f)-(v+f)\p_rf+\frac1rg^2.
\end{eqnarray*}
Since \eqref{axisymEuler} fulfills the null conditions of the potential equation near the light cone $\{r=t\}$, it follows from Lemma 3.1
and Remark 3.1 of Section 3
that the term $\p_r(v+f)$ and $v+f$ admit the  better decay in time near $\{r=t\}$. Note that both the quadratic error terms
in $Q_1$ and $Q_2$ contain the ``good" terms $\p_r(v+f)$ or $v+f$ (the left terms $\frac1rvf$ and $\frac1rg^2$
in $Q_1$ and $Q_2$ actually correspond to the fourth order error terms), and then the main parts of \eqref{axisymEuler}
are linear or linearly degenerate. Based on this basic observation, we can manage to
prove that \eqref{axisymEuler} will admit a global small data smooth solution $(v,f,g)$.
}

\vskip 0.2 true cm

Let's recall some remarkable works related to Theorem 1.1.
For the  2D or 3D compressible Euler equations of polytropic gases,
it is well known that smooth solution $(\rho, u)$ of \eqref{Euler} will generally blow up in finite time.
For examples,

{\bf $\bullet$} for a special class of initial data $(\rho_0(x), u_0(x))$,
T.~Sideris~\cite{Sideris85} has proved that the smooth solution $(\rho, u)$ of \eqref{Euler}  in three space dimensions can
develop singularities in finite time, and subsequently M.~Ramaha \cite{Rammaha89} has proved a corresponding
blowup result in two space dimensions;

{\bf $\bullet$} for the 2D compressible isentropic Euler equations, when the rotationally
invariant data are a perturbation of size $\ve>0$ of a rest state,  S.~Alinhac in \cite{Alinhac92} and \cite{Alinhac93}
establishes that the smooth solution blows up in finite time and  the lifespan
$T_{\ve}$ satisfies $\ds\lim_{\ve\to 0}\ve^2 T_{\ve}=\tau_{0}^2>0$;

{\bf $\bullet$} for the 2D and 3D  compressible Euler equations \eqref{Euler}, there are extensive
literature on the blowup or the blowup mechanism (including the formation of shocks) of small perturbed
smooth solution
$(\rho, u)$, one can see \cite{Alinhac99}, \cite{Christodoulou07}-\cite{CM}, \cite{DLY15}--\cite{Godin05}, \cite{HKSW}--\cite{John}, \cite{Lindblad08AJM}, \cite{LukSpeck16}--\cite{MS04}, \cite{Secchi02} and \cite{Sideris97}--\cite{Yin} (even for more general quasilinear wave equations).

\vskip 0.2 true cm
For the  2D or 3D compressible Euler equation of Chaplygin gases, so far there have been many
interesting results. For examples,

{\bf $\bullet$} when the  Chaplygin gases are isentropic and irrotational in two or three space dimensions, one can introduce a potential
function $\Phi(t,x)$ such that $u=\grad_x\Phi$ holds
and  \eqref{Euler} is written as a second order quasilinear wave equation of $\Phi$. In the 3D case,
this quasilinear wave equation
satisfies the null condition (its definition see \cite{Christodoulou86} and \cite{Klainerman}); in the 2D case, the equation
satisfies both the first and the second null condition (see their definitions in \cite{Alinhac01}).
Therefore, by the results in \cite{Christodoulou86}, \cite{Klainerman}  and \cite{Alinhac01}, \cite{LNS13},
we know that the small perturbed smooth solution $(\rho, u)$ of \eqref{Euler} exists globally when  3D or 2D  Chaplygin gases
are isentropic and irrotational.

{\bf $\bullet$} when the Chaplygin gases are spherically symmetric and non-isentropic, P.~Godin in
\cite{Godin07} prove that the global smooth symmetric solution exists for 3D non-isentropic compressible Euler equations.
Later on, the authors in \cite{DWY15} establish  a similar result in two space dimensions.
In addition, the authors in \cite{LeiWei17} show the global spherically symmetric solution
for 3D relativistic compressible Euler system of Chaplygin gases.

{\bf $\bullet$} when the 3D Chaplygin gases are isentropic and irrotational, the authors in \cite{LiYin18}
prove the global existence of small perturbed smooth solutions in the exterior domain with slip boundary condition $u\cdot\vec{n}=0$,
where $\vec{n}$ stands for the unit normal of the boundary. In fact, by proving
the variation version of global Keel-Smith-Sogge estimate, the authors in \cite{LiYin18}
have established more general result:  if 3D quasilinear wave equation satisfies the null condition,
then the  small data smooth solution exists globally  outside of a compact convex obstacle when the Neumann boundary condition
on the boundary is posed.

\vskip 0.1 true cm

We now give some comments on the proof of Theorem 1.1. When the 2D quasilinear wave equation satisfies
both the first and the second null conditions,  S.~Alinhac in  \cite{Alinhac01} have established the global existence of small data solution
by looking for a crucial ``ghost weight"  to derive a global energy estimate. Therefore, as a direct application of  \cite{Alinhac01},
when 2D Chaplygin gases are isentropic and irrotational, the small perturbed smooth solution $(\rho, u)$ of
\eqref{Euler} exists globally. If 2D Chaplygin gases are spherically symmetric and non-isentropic,
through looking for a new ``ghost weight'' and taking involved analysis,
the authors in \cite{DWY15} derive a global weighted energy estimate for
the full compressible  Euler system by utilizing both null conditions and the variable entropy.
As pointed out in \cite{DWY15}, since the full compressible Euler system can not be changed into a scalar quasilinear wave equation
because of the influence of variable entropy, the ghost weight introduced in \cite{Alinhac01} will not be applied directly
and thus a new  ghost weight should be chosen, meanwhile the transport equation on the entropy has to be treated
carefully by decomposing the whole energy into ``interior energy" and ``exterior energy". In the present paper, we shall focus on
the global solution problem of 2D compressible isentropic Euler equations of
Chaplygin gases when the rotationally
invariant data are a perturbation of size $\ve>0$ of a rest state.
Due to the strong effect of hydrodynamical waves on acoustical waves (i.e., the term $\frac{g^2}{r}$ in $Q_2$),
the time-decay of $v$ and $\p_rv$ away from the outgoing light cone (near $r\le t/2$ for large $t$) becomes worse,
see Lemma 2.5, Lemma 2.6 and Lemma 4.1 below. To overcome this essential difficulty,
we seek a suitable  transformation $v(t,r)=\vv(t,r)+G(t,r)$ so that the unknown function $\vv$ will
admit a better decay in time
away from the light cone meanwhile the function $G$ has some required ``good" properties in the process of deriving energy estimates
on $(\vv, f, g)$. For this purpose, we delicately choose $G$ to satisfy the nonlinear ODE:
\begin{equation}\label{G:def}
    (1+v)\p_rG+\frac{g^2}{r}=0
\end{equation}
with $G(t,\infty)=0$. Although \eqref{G:def} is nonlinear about the unknown function $G$ in terms of $v=\vv+G$, according to the estimates in Lemma~2.6, one knows that $G$ exists globally for $r\ge 0$ and any fixed $t\ge 0$.
In addition, in order to derive the global weighted energy estimate of  $(\vv, f, g)$,
it is required to search suitable vector fields and different weights near or away from the light cone. In fact, note that \eqref{axisymEuler}
are invariant under the scaling $(t,r)\rightarrow(\lambda t, \lambda r)$ for $\la\in\Bbb R$ and
the translation in time $t\rightarrow t+\tau$ for $\tau\in\Bbb R$, but \eqref{axisymEuler} are not invariant under radial translation.
Then this inspires that the vector fields $\Gamma\in\{\p_t,S=:t\p_t+r\p_r\}$
can be used  but the space derivative $\p_r$ can not be used directly in the energy estimates of $(\vv, f, g)$.
On the other hand, due to the requirements of argument techniques, near the light cone we will introduce the  auxiliary
weighted energy for $k\in\Bbb N$
\begin{equation}\label{energy0:def}
\begin{array}{l}
\displaystyle \cY_k(t)=:\sum_{|a|\le k-1}
\Big[\|<r-t>\chi_1(\frac{r}{<t>})\p_t\Gamma^av(t,r)\|_{L^2}+\|<r-t>\chi_1(\frac{r}{<t>})\p_r\Gamma^av(t,r)\|_{L^2}\\
\displaystyle \qquad\qquad+\|<r-t>\chi_1(\frac{r}{<t>})\p_t\Gamma^af(t,r)\|_{L^2}+\|<r-t>\chi_1(\frac{r}{<t>})\p_r\Gamma^af(t,r)\|_{L^2}\Big]
\end{array}
\end{equation}
together with some kind of ``ghost weight" and other suitable multiplier; meanwhile, away form the light cone we introduce the
auxiliary weighted energy
\begin{equation}\label{enerhy1:def}
\begin{array}{l}
\displaystyle \cX_k(t)=:<t>\sum_{|a|\le k-1}\Big[\|\chi_0(\frac{r}{<t>})\p_t\Gamma^a\vv(t,r)\|_{L^2}
+\|\chi_0(\frac{r}{<t>})\p_r\Gamma^a\vv(t,r)\|_{L^2}\\
\displaystyle \qquad\qquad+\|\chi_0(\frac{r}{<t>})\p_t\Gamma^af(t,r)\|_{L^2}+\|\chi_0(\frac{r}{<t>})(\p_r+\frac1r)\Gamma^af(t,r)\|_{L^2}\Big]\\
\displaystyle \qquad\qquad +<t>\sum_{|b|\le k-2}
\Big[\|\chi_0(\frac{r}{<t>})(\p_r+\frac1r)\p_r\Gamma^b\vv(t,r)\|_{L^2}\\
\displaystyle \qquad\qquad +\|\chi_0(\frac{r}{<t>})\p_r(\p_r+\frac1r)\Gamma^bf(t,r)\|_{L^2}\Big],
\end{array}
\end{equation}
where $\chi_0(s)$ and $\chi_1(s)$ are smooth cutoff functions satisfying $0\le\chi_0\le1$ with
\begin{equation}\label{cutoff:def}
\chi_0(s)+\chi_1(s)=1,\quad\chi_0(s)=\left\{
\begin{aligned}
    &1,\quad s\le\frac14,\\
    &0,\quad s\ge\frac12,
\end{aligned}
\right.
\end{equation}
for $s\in\Bbb R$.
Moreover, both $\frac{|\chi'_0|^2}{\chi_0}$ and $\frac{|\chi'_1|^2}{\chi_1}$ are bounded.
Here we emphasize that due to
$(\p_r+\frac1r)\Gamma^af(t,r)=div \Gamma^au$, $(\p_r+\frac1r)\p_r\Gamma^b\vv(t,r)=\Delta\Gamma^b\vv$ and $\p_r(\p_r+\frac1r)\Gamma^bf(t,r)=\grad\dive\Gamma^bu$, then each term in $\cX_k(t)$ is finite for
any fixed time $t\ge 0$ and smooth solution $(v, f, g)$ of \eqref{axisymEuler}.
In order to compensate the lack of the function $g$ and its regularities in $\cX_k(t)$, and to derive the $L^{\infty}$ estimate of $\p_rg$
by Sobolev imbedding theorem $W^{1,3}(\Bbb R^2)\subset L^{\infty}(\Bbb R^2)$ (see Lemma~4.2 for details), we shall introduce the  $W^{1,3}$ norm of the the specific vorticity
\begin{equation*}
    W_k(t)=:\sum_{|a|\le k}
    \big[\|\Gamma^aw(t,r)\|_{L^3}+\|\p_r\Gamma^aw(t,r)\|_{L^3}\big],
\end{equation*}
where $w(t,r)=\frac{\curl u(t,x)}{\rho(t,r)}=(1+v(t,r))(\p_r+\frac1r)g(t,r)$.
As in \cite{Speck16book}, it is easy to verify that $w(t,r)$ satisfies
\begin{equation}\label{curl:eqn}
    \p_tw+f\p_rw=0.
\end{equation}
Based on the preparations above and by involved analysis together with some delicate observations,
we eventually derive the global energy estimate of $(\vv, f, g)$ and further complete the proof of Theorem 1.1
by the continuation method.

This paper is organized as follows: In Section 2, several basic results including some Sobolev-type embedding
inequalities and Hardy-type inequalities are listed or derived. In Section 3, by the
null condition structures, the estimates of weighted $L^\infty$, $L^2$ norms and better decay property in time near the light
cone will be established for the smooth solution $(v, f, g)$ of \eqref{axisymEuler}. Away from the light cone, $L^\infty$ and $L^2$
time-decay of acoustical and hydrodynamical waves for problem \eqref{axisymEuler} are derived in Section 4.
In Section 5, by using S.~Alinhac's ``ghost weight" technique together with some other suitable multiplier, an elementary energy estimate
for the smooth solution of \eqref{axisymEuler} is obtained.
Thereafter, more careful energy estimates of $(v, f, g)$ near and away from the light cone will be given in Section 6 and 7, respectively.
In Section 8, the proof of Theorem~1.1 is eventually finished.

\vskip 0.1 true cm

Through the whole paper, we shall use the following convention:
\vskip 0.1 true cm
{\bf $\bullet$} $\nabla=(\p_{x_1}, \p_{x_2})$ for $x\in\Bbb R^2$.

\vskip 0.1 true cm
{\bf $\bullet$} $\Gamma\in\{\p_t,S=t\p_t+r\p_r\}$ and $\tilde\Gamma\in\{\p_t,\tilde S=:S+1\}$.

\vskip 0.1 true cm
{\bf $\bullet$} For the multi-indices $a,b\in\N_0^2=:\{0,1,2,\cdots\}^2$, then

$\Gamma^a=\p_t^{a_1}S^{a_2}$;

$a\le b$ means $a_1\le b_1$ and $a_2\le b_2$, while $a<b$ means $a_1<b_1$ or $a_2<b_2$ with $a\le b$;

$|a|=a_1+a_2$ and $a!=a_1!a_2!$;

Set $\sigma_{ab}=:\frac{(a+b)!}{a!b!}$.

\vskip 0.1 true cm
{\bf $\bullet$} The $L^p$ norm of $v(t,x)$ ($x\in \Bbb R^2$, $t\ge 0$) is
$\|v\|_{L^p}=\|v(t,x)\|_{L^p}=:(\int_{\Bbb R^2}|v(t,x)|^pdx)^{\f{1}{p}}$;

if $v(t,x)$ is symmetric, i.e., $v(t,x)=v(t,r)$, then $\|v\|_{L^p}=\|v(t,r)\|_{L^p(\Bbb R^2)}=(\int_0^\infty|v(t,r)|^prdr)^{\f{1}{p}}$.

\vskip 0.1 true cm

{\bf $\bullet$} For $k\in\N_0$ and the solution $(v,f,g)$ of \eqref{axisymEuler}, let
\begin{equation}\label{energy:def}
\begin{array}{l}
\displaystyle E_k(t)=:\sum_{|a|\le k}
\big[\|\Gamma^av(t,r)\|_{L^2}+\|\Gamma^af(t,r)\|_{L^2}+\|\Gamma^ag(t,r)\|_{L^2}\big]\\
\displaystyle \qquad\qquad+\sum_{|b|\le k-1}\|(\p_r+\frac1r)\Gamma^bg(t,r)\|_{L^2},
\end{array}
\end{equation}

where $(\p_r+\frac1r)\Gamma^bg(t,r)=curl \Gamma^bu$.

\vskip 0.2 true cm

\section{Some preliminaries}\label{sect2}

First, it follows from direct computation that for any smooth functions $\phi(t,r)$ and $\psi(t,r)$,
\begin{equation}
\begin{split}
&\tilde S(\phi\p_r\psi)=\phi\p_rS\psi+S\phi\p_r\psi,
\qquad \tilde S(\frac{\phi\psi}{r})=\frac1r\phi S\psi+\frac1rS\phi\psi,\no\\
&\tilde S\p_r\psi=\p_rS\psi, \qquad \tilde S(\frac\psi r)=\frac1rS\psi,
\qquad \tilde S\p_t\psi=\p_t S\psi.
\end{split}
\label{commutation}
\end{equation}
Through the paper, we always assume that for fixed integer $N\ge N_0=7$,
\begin{equation}
\begin{split}
\displaystyle E_N(t)+\cX_N(t)+\cY_N(t)\le M\ve(1+t)^{M'\ve},\\
\displaystyle E_{N-2}(t)+\cX_{N-2}(t)+\cY_{N-2}(t)+W_{N-3}(t)\le M\ve\le M\ve_0,
\end{split}
\label{assumption}
\end{equation}
where the large positive constants $M$ and $M'$ will be chosen later.
\vskip 0.1 true cm

{\bf Lemma 2.1.} {\it
For any function $\phi(t,x)$ which is supported in $\{x:|x|\le1+t\}$ for variable $x$ and for $t\ge 0$, the following Hardy-type inequality holds
\begin{equation}\label{hardyineq1}
    \Big\|\frac{\phi(t,x)}{<r-t>}\Big\|_{L^2}\ls\|\grad\phi(t,x)\|_{L^2},
\end{equation}
here and below $A\ls B$ means $A\le CB$ with generic positive constant $C$ which is independent of $t$, $\ve$ and $M,M'$.
}

{\bf Proof.} One can see the proof of \eqref{hardyineq1} in  \cite{Hans0}, here we omit it. \hfill $\square$
\vskip 0.2 true cm

{\bf Lemma 2.2.} {\it
For $a\in\N_0^2$,  if
$(v, f, g)$ is the solution of \eqref{axisymEuler}, then the following weighted Sobolev type inequalities hold
\begin{eqnarray}
<t>^\frac12<r-t>^\frac12|\chi_1\Gamma^af(t,r)|\ls E_{|a|}(t)+\cY_{|a|+1}(t),\label{f:pointwise}\\
<t>^\frac12<r-t>^\frac12|\chi_1\Gamma^av(t,r)|\ls E_{|a|}(t)+\cY_{|a|+1}(t),\label{v:pointwise}
\end{eqnarray}
where the definition of $\chi_1$ has been given in \eqref{cutoff:def}.}
\vskip 0.1 true cm
{\bf Remark 2.1.} {\it Set $cos\th =\f{x_1}{r}$ and $sin\th=\f{x_2}{r}$ with $\th\in [0, 2\pi]$. Then
$f(t,r)=u_1(r cos\th, r sin\th)cos\th+u_2(r cos\th, r sin\th)sin\th$ and $g(t,r)=u_2(r cos\th, r sin\th)cos\th-u_1(r cos\th, r sin\th)sin\th$.
Since $(\rho, u_1, u_2)(t,x)$ is a smooth solution of \eqref{EulerC1form} together with \eqref{ChaplyginGas},
one has $|\p_t^k\p_r^lf|+|\p_t^k\p_r^lg|\le C_{kl}(T)$ for $k,l\in\Bbb N_0$ and $t\in [0, T]$.}

\vskip 0.1 true cm

{\bf Proof.}
Note that $<t>+<r-t>\ls r$ on $\supp\chi_1(\frac{r}{<t>})$.
Since
\begin{eqnarray*}
  &&\left|<t><r-t>(\chi_1\Gamma^af)^2\right|=
  \left|\int_r^\infty<t>\frac{d}{dr'}(<r'-t>|\chi_1\Gamma^af|^2)dr'\right|\\
  &&\ls\int_r^\infty<r'-t>\chi_1|\chi'_1||\Gamma^af|^2dr'
  +\int_r^\infty <t><r'-t>\chi_1^2|\Gamma^af||\p_r\Gamma^af|dr'\\
  &&\qquad+\int_r^\infty <t>|\chi_1\Gamma^af|^2|\p_{r'}<r'-t>|dr',
\end{eqnarray*}
together with Cauchy-Schwartz inequality, the boundedness of $|\p_{r'}<r'-t>|$ and the definitions of $E_{|a|}(t)$, $\cY_{|a|+1}(t)$, this yields
\eqref{f:pointwise}. In addition, the proof of \eqref{v:pointwise} is similar, we then omit it. \hfill  $\square$

\vskip 0.1 true cm

{\bf Lemma 2.3.} {\it
For any smooth function $\phi(t,x)$ which has compact support for variable $x$, the following Sobolev inequalities hold
\begin{eqnarray}
    &\|\phi\|_{L^\infty}\ls\|\nabla\phi\|_{L^2}\ln^\frac12(2+t)+<t>^{-1}(\|\phi\|_{L^2}+\|\nabla^2\phi\|_{L^2}),\label{sobolevineq}\\
    &\|\phi\|_{L^\infty}\ls\|\phi\|_{L^3}+\|\nabla\phi\|_{L^3},\label{sobolev:W13}
\end{eqnarray}
where the definition of $\chi_0$ has been given in \eqref{cutoff:def}.}

\vskip 0.1 true cm

{\bf Proof.}
For the first inequality \eqref{sobolevineq}, one can see (3.4) of  \cite{Lei16}.
\eqref{sobolev:W13} comes from the Sobolev imbedding theorem $W^{1,3}(\Bbb R^2)\subset L^{\infty}(\Bbb R^2)$
directly.  \hfill $\square$

\vskip 0.1 true cm
{\bf Lemma 2.4.} {\it
For $C^1$ smooth function $\phi(t,r)$ with $\phi(t,0)=0$ and
admitting compact support for variable $r$,
the following Hardy type inequalities hold
\begin{eqnarray}
    &\|\frac\phi r\|_{L^2}\ls\|(\p_r+\frac1r)\phi\|_{L^2},\label{hardyineq2}\\
    &\|\frac1r\chi_0\phi\|_{L^2}\ls\|\chi_0(\p_r+\frac1r)\phi\|_{L^2}+<t>^{-1}\|\phi\|_{L^2}. \label{hardyineq3}
\end{eqnarray}
}

{\bf Proof.}
Direct computation yields
\begin{eqnarray*}
  \int_0^\infty|\frac1r\phi(t,r)|^2rdr=\int_0^\infty|r\phi(t,r)|^2d(\frac{-1}{2r^2})
  \le\left|\int_0^\infty\frac{\phi}{r}(\p_r+\frac1r)\phi rdr\right|\\
  \le\Big(\int_0^\infty|\frac1r\phi(t,r)|^2rdr\Big)^\frac12
  \Big(\int_0^\infty|(\p_r+\frac1r)\phi(t,r)|^2rdr\Big)^\frac12.
\end{eqnarray*}
Therefore, we have proved \eqref{hardyineq2}.

Next, we deal with \eqref{hardyineq3}.
By an analogous computation, we achieve
\begin{eqnarray*}
  &&\int_0^\infty|\frac1r\chi_0\phi(t,r)|^2rdr=\int_0^\infty|r\chi_0\phi(t,r)|^2d(\frac{-1}{2r^2})\\
  &&\ls\left|\int_0^\infty\chi_0^2\frac{\phi}{r}(\p_r+\frac1r)\phi rdr\right|
  +\left|\int_0^\infty<t>^{-1}\chi'_0\chi_0\phi^2dr\right|.
\end{eqnarray*}
Together with $\frac1C<t>\le r\le C<t>$ on $\supp\chi'_0$, this derives \eqref{hardyineq3}.
\hfill $\square$

\vskip 0.2 true cm

{\bf Lemma 2.5.} {\it
For $a\in\N_0^2$ and $|a|\le N-1$,  if
$(v, f, g)$ is  the solution of \eqref{axisymEuler}, then the following inequalities hold

(1) $|\Gamma^ag(t,r)|\ls E_{|a|+1}(t);$

(2) $|\frac1r\chi_0\Gamma^af(t,r)|\ls\|\chi_0(\p_r+\frac1r)\Gamma^af(t,r)\|_{L^\infty}
+<t>^{-1}\|\chi^\frac12_0\Gamma^af(t,r)\|_{L^\infty};$

(3) $<t>|\chi_0\Gamma^af(t,r)|\ls E_{|a|}(t)+\cX_{|a|+1}(t);$

(4) $<t>^{0.99}|\chi_0\Gamma^a\vv(t,r)|\ls E_{|a|}(t)+\cX_{|a|+2}(t)
+\cY_{|a|+1}(t)+\|\Gamma^aG(t,r)\|_{L^2}+\|\chi_1\p_r\Gamma^aG(t,r)\|_{L^2},$
where $\vv=v-G$, and the definition of $G$ is given in \eqref{G:def}.

Furthermore, for any fixed $t\ge0$, it holds that

(5) $\Gamma^ag(t,0)=\Gamma^af(t,0)=0$ with $|a|\le N-1$ and  $\p_r\Gamma^{a'}\vv(t,0)=0$ with $|a'|\le N-2$.}

\vskip 0.1 true cm
{\bf Remark 2.2.} {\it In view of the result (5), Lemma~2.4 can be applied to the terms
$\Gamma^ag(t,r), \Gamma^af(t,r)$ and $\p_r\Gamma^{a'}\vv(t,r)$ for $|a|\le N-1$ and $|a'|\le N-2$.}

\vskip 0.1 true cm

{\bf Proof.} At first, we prove (1). Due to
\begin{eqnarray*}
  &&|\Gamma^ag(t,r)|=\frac1r\left|\int_0^r\p_r(r'\Gamma^ag)dr'\right|\\
  &&=\frac1r\left|\int_0^r(\p_r+\frac{1}{r'})\Gamma^agr'dr'\right|\ls
  \frac1r\Big(\int_0^rr'dr'\Big)^\frac12\Big(\int_0^r|(\p_r+\frac{1}{r'})\Gamma^ag|^2r'dr'\Big)^\frac12,
\end{eqnarray*}
then (1) is proved.

In addition, the boundedness of $\int_0^\infty|(\p_r+\frac1r)\Gamma^ag|^2rdr$ implies that $\int_0^r|(\p_r+\frac{1}{r'})\Gamma^ag|^2r'dr'\rightarrow0$ as $r\rightarrow0$.
Then we find that $\Gamma^ag(t,0)=0$.
With the same analysis, we also have $\Gamma^af(t,0)=\p_r\Gamma^{a'}\vv(t,0)=0$.
Therefore, (5) is achieved.

Next, we deal with (2).
Similar to the proof of (1), we see that
\begin{equation}\label{f:sobolev}
\begin{array}{l}
\displaystyle |\chi_0r\Gamma^af(t,r)|=\left|\int_0^r\p_r(\chi_0r'\Gamma^af)dr'\right|\\
\displaystyle \ls\left|\int_0^r\chi_0(\p_r+\frac{1}{r'})\Gamma^afr'dr'\right|
+\left|\int_0^r<t>^{-1}\chi'_0\Gamma^afr'dr'\right|.
\end{array}
\end{equation}
Therefore,
\begin{eqnarray*}
  |\frac1r\chi_0\Gamma^af(t,r)|\ls\frac{1}{r^2}\int_0^rr'dr'
  \Big[\|\chi_0(\p_r+\frac1r)\Gamma^af(t)\|_{L^\infty}+<t>^{-1}\|\chi^\frac12_0\Gamma^af(t)\|_{L^\infty}\Big],
\end{eqnarray*}
which implies (2).

Applying the Cauchy-Schwartz inequality to \eqref{f:sobolev} yields
\begin{eqnarray*}
  |\chi_0\Gamma^af(t,r)|\ls\frac1r\Big(\int_0^rr'dr'\Big)^\frac12
  \Big[\|\chi_0(\p_r+\frac1r)\Gamma^af\|_{L^2}+<t>^{-1}\|\Gamma^af(t)\|_{L^2}\Big].
\end{eqnarray*}
Then we have proved (3).

Finally, we turn our attention to the proof of (4).
Let $\phi=\chi_0\Gamma^a\vv(t,r)$ in \eqref{sobolevineq}, one has that
\begin{equation}\label{v:sobolev}
    |\chi_0\Gamma^a\vv(t,r)|\ls\ln^\frac12(2+t)\Big[<t>^{-1}\|\Gamma^a\vv\|_{L^2}
    +\|\p_r\Gamma^a\vv\|_{L^2}+\|(\p_r+\frac1r)\p_r(\chi_0\Gamma^a\vv)\|_{L^2}\Big].
\end{equation}
We now deal with the last term in \eqref{v:sobolev}. In fact, it is easy to find that
\begin{equation}\label{d-0}
\begin{array}{l}
    \|(\p_r+\frac1r)\p_r(\chi_0\Gamma^a\vv)\|_{L^2}\ls\|\chi_0(\p_r+\frac1r)\p_r\Gamma^a\vv\|_{L^2}
    +\|(<t>^{-2}\chi''_0+<t>^{-1}\chi'_0\frac1r)\Gamma^a\vv\|_{L^2}\\
   \qquad\qquad  +\|(<t>^{-1}|\chi'_0|+\chi_0\frac1r)\p_r\Gamma^a\vv\|_{L^2}
    +\|<t>^{-1}\chi'_0(\p_r+\frac1r)\Gamma^a\vv\|_{L^2}.
\end{array}
\end{equation}
Applying \eqref{hardyineq3} to $\|\chi_0\frac1r\p_r\Gamma^a\vv\|_{L^2}$, we deduce that
\begin{equation}\label{d-1}
\begin{array}{l}
    \|(\p_r+\frac1r)\p_r(\chi_0\Gamma^a\vv)\|_{L^2}\ls\|\chi_0(\p_r+\frac1r)\p_r\Gamma^a\vv\|_{L^2}
    +<t>^{-1}\|\p_r\Gamma^a\vv\|_{L^2}\\
   \qquad\qquad   +<t>^{-2}\|\Gamma^a\vv(|\chi'_0|+|\chi''_0|)\|_{L^2}.
\end{array}
\end{equation}
Substituting \eqref{d-0} and \eqref{d-1} into \eqref{v:sobolev} yields (4). \hfill $\square$

\vskip 0.2 true cm

Next we estimate $G$ in \eqref{G:def}.

\vskip 0.1 true cm

{\bf Lemma 2.6.} {\it
For $a, a'\in\N_0^2$ with $|a|\le N-2$ and $|a'|\le N$, the following inequalities hold

(1) $|\Gamma^aG|\ls\int_0^\infty|\p_r\Gamma^aG|dr\ls E_{|a|+1}(t);$

(2) $|r\p_r\Gamma^aG(t,r)|\ls E_{|a|+1}(t);$

(3) $\|\Gamma^{a'}G\|_{L^2}+\|r\p_r\Gamma^{a'}G\|_{L^2}\ls E_{|a'|}(t).$
}
\vskip 0.1 true cm

{\bf Proof.} (1) From the third equation of \eqref{axisymEuler}, we have
\begin{equation*}
    \p_t(rg)+f\p_r(rg)=0.
\end{equation*}
Define the characteristic curve $R(t,r)$ starting from $(0,r)$ as follows:
\begin{equation*}
    \frac{d}{dt}R(t,r)=f(t,R(t,r)),\quad R(0,r)=r.
\end{equation*}
Then $R(t,r)g(t,R(t,r))=rg(0,r)$ holds.
In view of \eqref{f:pointwise} and Lemma 2.5,
we deduce that $|f(t,R(t,r))|$ $\ls<t>^{-\frac12}[E_0(t)+\cX_1(t)+\cY_1(t)]$.
Then it is easy to get $|R(t,r)-r|\le CM\ve_0<t>^\frac12$.
By the assumption $\supp g_0(r)\subseteq\{r\le1/8\}$,
one has $\supp g(t,r)\subseteq\{(t,r):r\le<t>/4\}$.
Consequently, from the definition of $G$, we get
\begin{equation}\label{G:integral}
    G(t,r)=\int_r^\infty\frac{g^2}{1+v}(t,r')\frac{dr'}{r'}.
\end{equation}
This implies
\begin{equation}\label{g:support}
\begin{array}{l}
\displaystyle \supp G(t,r)\subseteq\{(t,r):r\le<t>/4\},\\
\displaystyle \chi_0G\equiv G,\quad\chi_1G\equiv0,\quad\chi_0g\equiv g,\quad\chi_1g\equiv0.
\end{array}
\end{equation}
Hence, for any fixed $t\ge0$ we have $\ds\lim_{r\rightarrow\infty}\Gamma^aG(t,r)=0$.
Applying the commutation identities \eqref{commutation},
we get the following  equation for the higher order derivatives of $G$
\begin{equation}\label{G:eqn1}
    (1+v)\p_r\Gamma^aG+\sum_{b+c=a,c<a}\sigma_{bc}(\Gamma^b\vv+\Gamma^bG)\p_r\Gamma^cG
    +\frac1r\sum_{b+c=a}\sigma_{bc}\Gamma^bg\Gamma^cg=0.
\end{equation}
Therefore,
\begin{equation}\label{G:L1norm1}
\begin{array}{l}
\displaystyle (1-\|\chi_0\vv\|_{L^\infty}-\|G\|_{L^\infty})\int_0^\infty|\p_r\Gamma^aG|dr\\
\displaystyle\ls\sum_{b+c=a,c<a}(\|\chi_0\Gamma^b\vv\|_{L^\infty}+\|\Gamma^bG\|_{L^\infty})
\int_0^\infty|\p_r\Gamma^cG|dr+\sum_{b+c=a}\int_0^\infty|\Gamma^bg\Gamma^cg|\frac1rdr
\end{array}
\end{equation}
for $|a|\neq0$, and
\begin{equation}\label{G:L1norm1:a=0}
  (1-\|\chi_0\vv\|_{L^\infty}-\|G\|_{L^\infty})\int_0^\infty|\p_rG|dr
  \ls\int_0^\infty\frac{g^2}{r}dr.
\end{equation}
With the help of the Cauchy-Schwartz inequality and Hardy inequality \eqref{hardyineq2},
one has
\begin{equation}\label{g:hardy1}
\begin{array}{l}
\displaystyle \int_0^\infty|\Gamma^bg\Gamma^cg|\frac1rdr\ls
\int_0^\infty\frac{|\Gamma^bg|^2+|\Gamma^cg|^2}{r}dr\\
\displaystyle \ls\int_0^\infty|(\p_r+\frac1r)\Gamma^bg|^2rdr
+\int_0^\infty|(\p_r+\frac1r)\Gamma^cg|^2rdr.
\end{array}
\end{equation}
Collecting Lemma~2.5 and \eqref{G:L1norm1}-\eqref{g:hardy1},
we arrive at
\begin{equation}\label{G:L1norm2}
\begin{array}{l}
\displaystyle \int_0^\infty|\p_r\Gamma^aG|dr\\
\displaystyle \ls\sum_{b+c=a,c<a}(M\ve_0+\int_0^\infty|\p_r\Gamma^bG|dr)\int_0^\infty|\p_r\Gamma^cG|dr
+E_{|a|+1}(t).
\end{array}
\end{equation}

Next we apply the method of induction on the multi-index $a$ to show (1). In fact, for $|a|\le N-4$, one easily gets (1).
When $|a|\le N-2$, due to $|b|+2\le N-2$ or $|c|+2\le N-2$, then one has $E_{|b|+1}(t)+\cX_{|b|+2}(t)\le CM\ve_0$ or
$E_{|c|+1}(t)+\cX_{|c|+2}(t)\le CM\ve_0$.
Therefore, \eqref{G:L1norm2} always yields
\begin{equation*}
    \int_0^\infty|\p_r\Gamma^aG|dr\ls M\ve_0\sum_{c\le a}
    \int_0^\infty|\p_r\Gamma^cG|dr+E_{|a|+1}(t),
\end{equation*}
which derives (1) for sufficiently small $\ve_0>0$.

(2) By an analogous or even easier analysis as in (1), one can obtain (2).

(3) We conclude that
\begin{equation*}
    \|\Gamma^{a'}G\|^2_{L^2}=\int_0^\infty|\Gamma^{a'}G|^2d(\frac{r^2}{2})=
    \left|\int\Gamma^{a'}G\p_r\Gamma^cGr^2dr\right|\le\|\Gamma^{a'}G\|_{L^2}\|r\p_r\Gamma^{a'}G\|_{L^2}.
\end{equation*}
Note that $[\frac{|a'|}{2}]+2\le[\frac{|N-1|}{2}]+2\le N-2$
for $N\ge N_0$.
Then applying the method of induction to \eqref{G:eqn1}, we easily find that
\begin{equation*}
    \|\Gamma^{a'}G\|_{L^2}\le\|r\p_r\Gamma^{a'}G\|_{L^2}\ls E_{|a'|}(t)\big[E_{[\frac{|a'|}{2}]+1}(t)+\cX_{[\frac{|a'|}{2}]+2}(t)+\cY_{[\frac{|a'|}{2}]+1}(t)\big]\ls E_{|a'|}(t),
\end{equation*}
which implies (3). \hfill $\square$

\vskip 0.2 true cm

\section{On the analysis   of solutions near the light cone}\label{sect3}

In this section, for the smooth solution $(v, f, g)$ of \eqref{axisymEuler}, we will make full use of the inherent
null condition structure near the light cone
to establish some time-decay estimates.
\vskip 0.1 true cm

{\bf Lemma 3.1.} {\it For $a, a'\in\N_0^2$ with $|a|\le N-2$ and $|a'|\le N-1$, if
$(v, f, g)$ is the solution of \eqref{axisymEuler}, then
the following estimates of $L^\infty$ and $L^2$ norms for $v+f$ hold
\begin{align}
  &<t>^\frac32|\chi_1\p_r\Gamma^a(v+f)|\ls E_{|a|+1}(t)+\cY_{|a|+2}(t),\label{goodLinfty}\\
  &<t>\|\chi_1\p_r\Gamma^{a'}(v+f)\|_{L^2}\ls E_{|a'|+1}(t)\label{goodL2norm}
\end{align}
and
\begin{equation}\label{goodLinfty'}
\begin{split}
 & <r-t>^{-1}\chi_1|\Gamma^a(v+f)|\ls\|\chi_1\p_r\Gamma^a(v+f)\|_{L^\infty}
  +<t>^{-1}(\|\chi_1\Gamma^a(v+f)\|_{L^\infty}\no\\
 & \qquad \qquad\qquad \qquad \qquad \qquad +\|\chi_0\Gamma^a(\vv+f)\|_{L^\infty}).
\end{split}
\end{equation}
}

\vskip 0.1 true cm
{\bf Remark 3.1.} {\it  Near the light cone, since the 2D Chaplygin gases will be irrotational,
we can introduce the potential function $\Phi$ so that $u=\na\Phi$ for \eqref{EulerC1form}.
In this case, \eqref{EulerC1form} can be changed as the following quasilinear wave equation:
\begin{equation}\label{potential:eqn}
\begin{split}
&\p_t^2\Phi+2\ds\sum_{i=1}^2\p_i\Phi\p_t\p_i\Phi+\ds\sum_{i,j=1}^2\p_i\Phi\p_j\Phi\p_{ij}^2\Phi
-(1+2\p_t\Phi+|\nabla\Phi|^2)\Delta\Phi=0.
\end{split}
\end{equation}
It is easy to verify that both the first and second null condition hold for \eqref{potential:eqn}.
By the results in \cite{Alinhac01}, the $L^\infty$ and $L^2$ norms  of ``good derivative" $(\p_t+\p_r)\Phi$
near the light cone admit better time-decays. On the other hand, with the help of $f=\p_r\Phi$,
integrating the second equation of \eqref{axisymEuler} with respect to $r$ yields
$v=\p_t\Phi+l.o.t.$ near the light cone.
Here $l.o.t.$ means the lower order terms that decay faster in time.
Therefore, $v+f=(\p_t+\p_r)\Phi+l.o.t.$ should have better time-decay properties near the light cone
as indicated in Lemma 3.1.
}

\vskip 0.1 true cm

{\bf Proof.}
Acting $\Gamma^a$ on two sides of the first two equations in \eqref{axisymEuler} and applying \eqref{commutation}, we have
\begin{equation}\label{highordereqn}
\left\{
\begin{aligned}
    &\p_t\Gamma^av-(\p_r\Gamma^af+\frac1r\Gamma^af)=\tilde\Gamma^aQ_1,\\
    &\p_t\Gamma^af-\p_r\Gamma^av=\tilde\Gamma^aQ_2.
\end{aligned}
\right.
\end{equation}

Direct computation shows that
\begin{equation}\label{good:idetity}
\begin{split}
    &(t+r)\p_r\Gamma^a(v+f)=t\p_r\Gamma^a(v+f)+S\Gamma^a(v+f)-t\p_t\Gamma^a(v+f)\\
    &=-t\tilde\Gamma^a(Q_1+Q_2)+S\Gamma^av+S\Gamma^af-\frac tr\Gamma^af,
\end{split}
\end{equation}
where
\begin{equation*}
    \chi_1\tilde\Gamma^a(Q_1+Q_2)=\chi_1\sum_{b+c=a}\sigma_{bc}
    [\Gamma^b(v-f)\p_r\Gamma^c(v+f)+\frac1r\Gamma^bv\Gamma^cf].
\end{equation*}

Since $g\equiv0$ holds on $\supp\chi_1$, which has been proved in Lemma 2.6,
this leads to that $\chi_1\tilde\Gamma^aQ_2$ does not contain the function $g$.
By $\chi_1=\chi_1(\chi_0+\chi_1)$, we arrive at
\begin{equation}\label{good:ineq1}
\begin{split}
    \sum_{b+c=a}&|\chi_1\frac1r\Gamma^bv\Gamma^cf|\le\sum_{b+c=a,|b|\le|c|}
    |\frac1r(\chi_0\Gamma^b\vv+\chi_1\Gamma^bv)\chi_1\Gamma^cf|
    +\sum_{b+c=a,|c|\le|b|}|\frac1r\chi_1\Gamma^bv(\chi_0\Gamma^cf+\chi_1\Gamma^cf)|\\
    &\ls<t>^{-\frac32}\big[E_{|a|}(t)+\cY_{|a|+1}(t)\big]
    \big[E_{[\frac{|a|}{2}]}(t)+\cX_{[\frac{|a|}{2}]+2}(t)+\cY_{[\frac{|a|}{2}]+1}(t)\big],
\end{split}
\end{equation}
where $[s]=\sup\{i\in\N:i\le s\}$. In addition, for $|a|\le N-2$ and $N\ge N_0$, it is easy to check that $[\frac{|a|}{2}]+2\le N-2$.

Analogously, we find that
\begin{equation}\label{good:ineq2}
\begin{split}
    &\sum_{b+c=a}|\chi_1\Gamma^b(v-f)\p_r\Gamma^c(v+f)|\ls\sum_{b+c=a}
    \left|\chi_1\p_r\Gamma^c(v+f)\big[\chi_0\Gamma^b(\vv-f)+\chi_1\Gamma^b(v-f)\big]\right|\\
    &\qquad \ls<t>^{-\frac12}\sum_{b+c=a}|\chi_1\p_r\Gamma^c(v+f)|
    \big[E_{|b|}(t)+\cX_{|b|+2}(t)+\cY_{|b|+1}(t)\big].
\end{split}
\end{equation}
Then substituting \eqref{good:ineq1} and \eqref{good:ineq2} into \eqref{good:idetity} yields \eqref{goodLinfty}.

The proof of \eqref{goodL2norm} is similar. Indeed, by $N\ge N_0$ and $|b|+|c|\le|a'|\le N-1$,
one has $|b|\le N-2$ or $|c|\le N-2$. Therefore, we see that
\begin{eqnarray*}
    \|\chi_1\Gamma^b(v-f)\p_r\Gamma^c(v+f)\|_{L^2}\ls\sum_{b+c=a',|c|\le N-2}
    \|\chi_1\p_r\Gamma^c(v+f)\|_{L^\infty}\|\Gamma^b(v-f)\|_{L^2}\\
    +<t>^{-\frac12}\sum_{b+c=a',|b|\le N-2}\|\chi_1\p_r\Gamma^c(v+f)\|_{L^2}
    \big[E_{|b|}(t)+\cX_{|b|+2}(t)+\cY_{|b|+1}(t)\big].
\end{eqnarray*}
This, together with \eqref{goodLinfty}, yields \eqref{goodL2norm}.

By the finite propagation speed property of hyperbolic system,
one easily knows that $\Gamma^a(v+f)$ is supported in $\{(t,r):r\le1+t\}$.
Therefore,
\begin{equation}\label{goodLinfty''}
\begin{array}{l}
    <r-t>^{-1}\chi_1|\Gamma^a(v+f)|=<r-t>^{-1}
    \left|\int_r^{1+t}\p_r(\chi_1\Gamma^a(v+f))dr'\right|\\ \ls<t>^{-1}\|\chi'_1\Gamma^a(v+f)\|_{L^\infty}
    +\|\chi_1\p_r\Gamma^a(v+f)\|_{L^\infty}.
\end{array}
\end{equation}

Since $\frac{|\chi'_1|^2}{\chi_1}$ is bounded and $v=\vv$ holds on $\supp\chi_1$, we arrive at
\begin{eqnarray*}
    \|\chi'_1\Gamma^a(v+f)\|_{L^\infty}\ls\|\chi_1^\frac12\Gamma^a(v+f)\|_{L^\infty}
    \ls\|\chi_1^\frac32\Gamma^a(v+f)\|_{L^\infty}+\|\chi_1^\frac12\chi_0\Gamma^a(v+f)\|_{L^\infty}\\
    \ls\|\chi_1\Gamma^a(v+f)\|_{L^\infty}+\|\chi_1^\frac12\chi_0\Gamma^a(\vv+f)\|_{L^\infty}
    \ls\|\chi_1\Gamma^a(v+f)\|_{L^\infty}+\|\chi_0\Gamma^a(\vv+f)\|_{L^\infty}.
\end{eqnarray*}
Substituting this into \eqref{goodLinfty''} yields \eqref{goodLinfty'}. \hfill $\square$

\vskip 0.1 true cm

{\bf Lemma 3.2.} {\it For $a\in\N_0^2$ with $|a|\le N-3$, if
$(v, f, g)$ is  the  solution of \eqref{axisymEuler}, then for small $\ve_0>0$,
the following estimates of weighted $L^\infty$ norms hold
\begin{equation}\label{nearLC:f:pw}
\begin{split}
&|<r-t>^\frac32\chi_1\p_r\Gamma^af(t,r)|+|<r-t>^\frac32\chi_1\p_t\Gamma^af(t,r)|\\
&\qquad \ls<t>^{-\frac12}\big[E_{|a|+1}(t)+\cY_{|a|+2}(t)\big]
\end{split}
\end{equation}
and
\begin{equation}\label{nearLC:v:pw}
\begin{split}
&|<r-t>^\frac32\chi_1\p_r\Gamma^av(t,r)|+|<r-t>^\frac32\chi_1\p_t\Gamma^av(t,r)|\\
&\ls<t>^{-\frac12}\big[E_{|a|+1}(t)+\cY_{|a|+2}(t)\big].
\end{split}
\end{equation}
}

\vskip 0.1 true cm

{\bf Proof.}
Direct computation yields
\begin{equation}\label{weighteqn1}
\begin{array}{l}
\displaystyle (t^2-r^2)\p_r\Gamma^av=t^2(\p_t\Gamma^af-\tilde\Gamma^aQ_2)
-rS\Gamma^av+tr\p_t\Gamma^av\\
\displaystyle =tS\Gamma^af-tr\p_r\Gamma^af-t^2\tilde\Gamma^aQ_2
-rS\Gamma^av+tr\p_t\Gamma^av\\
\displaystyle =tr\tilde\Gamma^aQ_1-t^2\tilde\Gamma^aQ_2+tS\Gamma^af-rS\Gamma^av+t\Gamma^af.
\end{array}
\end{equation}
Similarly, we also have
\begin{equation}\label{weighteqn1'}
\begin{array}{l}
\displaystyle (t^2-r^2)\p_t\Gamma^av=-r^2\tilde\Gamma^aQ_1+tr\tilde\Gamma^aQ_2
+tS\Gamma^av-rS\Gamma^af-r\Gamma^af,\\
\displaystyle (t^2-r^2)\p_t\Gamma^af=tr\tilde\Gamma^aQ_1-r^2\tilde\Gamma^aQ_2
-rS\Gamma^av+tS\Gamma^af+t\Gamma^af,\\
\displaystyle (t^2-r^2)\p_r\Gamma^af=-t^2\tilde\Gamma^aQ_1+tr\tilde\Gamma^aQ_2
+tS\Gamma^av-rS\Gamma^af+\frac{t^2}{r}\Gamma^af.
\end{array}
\end{equation}
By Leibniz's formula and the commutation identities \eqref{commutation}, one has
\begin{equation}\label{highorderQ1}
    \tilde\Gamma^aQ_1=\sum_{b+c=a}\sigma_{bc}[\Gamma^bv\p_r\Gamma^c(v+f)
    -\Gamma^b(v+f)\p_r\Gamma^cv+\frac1r\Gamma^bv\Gamma^cf]
\end{equation}
and
\begin{equation}\label{highorderQ2}
    \chi_1\tilde\Gamma^aQ_2=\chi_1\sum_{b+c=a}\sigma_{bc}
    [\Gamma^bv\p_r\Gamma^c(v+f)-\Gamma^b(v+f)\p_r\Gamma^cf].
\end{equation}
By an analogous analysis as in Lemma~3.1, in order to prove \eqref{nearLC:f:pw}
and \eqref{nearLC:v:pw},
it suffices only to deal with the term $\Gamma^b(v+f)\p_r\Gamma^cv$  in \eqref{highorderQ1}.
Note that
\begin{equation*}
\begin{array}{l}
    <t>|<r-t>^\frac12\chi_1\Gamma^b(v+f)\p_r\Gamma^cv|\\
    \ls <t>|<r-t>^\frac12\chi_1\Gamma^b(v+f)\chi_1\p_r\Gamma^cv|
   +<t>^\frac32|\chi_1\Gamma^b(v+f)\chi_0\p_r\Gamma^c\vv|\\
    \ls<t>^{-\frac12}|<r-t>^\frac32\chi_1\p_r\Gamma^cv|\big[E_{|b|+1}(t)
    +\cX_{|b|+2}(t)+\cY_{|b|+2}(t)\big]\\
    \quad +<t>^{-0.99}\big[E_{|b|}(t)+\cY_{|b|+1}(t)\big]
\big[E_{|c|+2}(t)+\cX_{|c|+3}(t)\big],
\end{array}
\end{equation*}
which implies
\begin{equation}\label{Q1:pw:ineq}
\begin{split}
<t>|<r-t>^\frac12\chi_1\tilde\Gamma^aQ_1|&\ls M\ve_0\sum_{c\le a}
    |<r-t>^\frac32\chi_1\p_r\Gamma^cv|\\
&+<t>^{-\frac12}\big[E_{|a|}(t)+\cY_{|a|+1}(t)\big].
\end{split}
\end{equation}
In addition, for $\chi_1\tilde\Gamma^aQ_2$, the same inequality as \eqref{Q1:pw:ineq} can be easily obtained.
Therefore, when $\ve_0>0$ is small,
we achieve \eqref{nearLC:f:pw}--\eqref{nearLC:v:pw}. \hfill $\square$

\vskip 0.1 true cm

{\bf Lemma 3.3.} {\it
For any integer $k\in\N$ with $k\le N$, if
$(v, f, g)$ is  the  solution of \eqref{axisymEuler}, then for small $\ve_0>0$,
the following auxiliary weighted energy inequality holds
\begin{equation}\label{WT:energy1}
    \cY_k(t)\ls E_k(t).
\end{equation}
}

\vskip 0.1 true cm

{\bf Proof.} Set
\begin{eqnarray*}
    \cE_1^a(t)=:\|<r-t>\chi_1\p_t\Gamma^av(t,r)\|_{L^2}+\|<r-t>\chi_1\p_r\Gamma^av(t,r)\|_{L^2}\\
    +\|<r-t>\chi_1\p_t\Gamma^af(t,r)\|_{L^2}+\|<r-t>\chi_1\p_r\Gamma^af(t,r)\|_{L^2}.
\end{eqnarray*}
In view of \eqref{weighteqn1} and \eqref{weighteqn1'}, it is not hard to see that
\begin{equation}\label{WE11}
    \cE_1^a(t)\ls E_{|a|+1}(t)+<t>\|\chi_1\tilde\Gamma^aQ_1\|_{L^2}
    +<t>\|\chi_1\tilde\Gamma^aQ_2\|_{L^2}.
\end{equation}
Since the treatments on $\chi_1\tilde\Gamma^aQ_1$ and $\chi_1\tilde\Gamma^aQ_2$ are analogous,
it suffices only to treat $\chi_1\tilde\Gamma^aQ_1$.

First, we deal with the term $\Gamma^bv\p_r\Gamma^c(v+f)$ in $\tilde\Gamma^aQ_1$.
Due to $|b|+|c|=|a|\le N$ with $N\ge N_0$,
one has $|b|\le N-2$ or $|c|\le N-2$.
Therefore, we achieve
\begin{equation}\label{WE12}
\begin{array}{l}
\displaystyle \sum_{b+c=a}\|\chi_1\Gamma^bv\p_r\Gamma^c(v+f)\|_{L^2}\ls
\sum_{b+c=a,|c|\le N-2}\|\chi_1\p_r\Gamma^c(v+f)\|_{L^\infty }\|\Gamma^bv\|_{L^2}\\
\displaystyle \qquad+\sum_{b+c=a,|b|\le N-2}\|\chi_1\p_r\Gamma^c(v+f)\|_{L^2}
\big[\|\chi_0\Gamma^b\vv\|_{L^\infty}+\|\chi_1\Gamma^bv\|_{L^\infty}\big]\\
\displaystyle \ls<t>^{-\frac32}\sum_{b+c=a,|c|\le N-2}E_{|b|}(t)\big[E_{|c|+1}(t)
+\cY_{|c|+2}(t)\big]\\
\displaystyle \qquad+<t>^{-\frac32}\sum_{b+c=a,|b|\le N-2}E_{|c|+1}(t)\big[E_{|b|}(t)
+\cX_{|b|+2}(t)+\cY_{|b|+1}(t)\big].
\end{array}
\end{equation}

Analogously for the term $\Gamma^c(v+f)\p_r\Gamma^bv$ in $\tilde\Gamma^aQ_1$, we have
\begin{equation}\label{WE13}
\begin{split}
&\sum_{b+c=a}\|\chi_1\Gamma^c(v+f)\p_r\Gamma^bv\|_{L^2}\\
&=\sum_{b+c=a,|c|\le N-2}\|<r-t>\chi_1\p_r\Gamma^bv\|_{L^2}
\|<r-t>^{-1}\chi_0\Gamma^c(\vv+f)+<r-t>^{-1}\chi_1\Gamma^c(v+f)\|_{L^\infty}\\
&\quad+\sum_{b+c=a,|b|\le N-2}\|<r-t>\chi_1\p_r\Gamma^bv\|_{L^\infty}
\|<r-t>^{-1}\chi_0\Gamma^c(v+f)+<r-t>^{-1}\chi_1\Gamma^c(v+f)\big]\|_{L^2}\\
&\ls<t>^{-\frac32}\sum_{b+c=a,|c|\le N-2}
\cE_1^b(t)\big[E_{|c|+1}(t)+\cX_{|c|+2}+\cY_{|c|+2}\big]\\
&\quad+<t>^{-\frac32}\sum_{b+c=a,|b|\le N-2}\big[E_{|b|+1}(t)+\cY_{|b|+2}(t)\big]
\big[E_{|c|}(t)+<t>\Big\|\chi_1\frac{\Gamma^c(v+f)}{<r-t>}\Big\|_{L^2}\big].
\end{split}
\end{equation}

Applying Hardy inequality \eqref{hardyineq1} to the last term in \eqref{WE13}, we obtain
\begin{equation*}
    \Big\|\chi_1\frac{\Gamma^c(v+f)}{<r-t>}\Big\|_{L^2}\ls
    \|\chi_1\p_r\Gamma^c(v+f)\|_{L^2}+<t>^{-1}\|\chi'_1\Gamma^c(v+f)\|_{L^2}.
\end{equation*}

This, together with (3) in Lemma~2.6 and \eqref{WE13},  yields
\begin{equation}\label{WE14}
\begin{array}{l}
\displaystyle \sum_{b+c=a}\|\chi_1\Gamma^c(v+f)\p_r\Gamma^bv\|_{L^2}\\
\displaystyle \ls<t>^{-\frac32}\sum_{b+c=a,|c|\le N-2}
\cE_1^b(t)\big[E_{|c|+1}(t)+\cX_{|c|+2}+\cY_{|c|+2}\big]\\
\displaystyle \qquad+<t>^{-\frac32}\sum_{b+c=a,|b|\le N-2}
E_{|c|+1}(t)\big[E_{|b|+1}(t)+\cY_{|b|+2}(t)\big].
\end{array}
\end{equation}

Finally, for the term $\frac1r\Gamma^bv\Gamma^cf$ in $\tilde\Gamma^aQ_1$, it is easy to get
\begin{equation*}
    \sum_{b+c=a}\|\chi_1\frac1r\Gamma^bv\Gamma^cf\|_{L^2}\ls
    <t>^{-\frac32}E_{|a|}\big[E_{[\frac{|a|}{2}]}+\cY_{[\frac{|a|}{2}]+1}\big].
\end{equation*}

Thus, collecting \eqref{WE11} with \eqref{WE12} and \eqref{WE14} leads to
\begin{equation*}
    <t>\|\chi_1\tilde\Gamma^aQ_1\|_{L^2}+<t>\|\chi_1\tilde\Gamma^aQ_2\|_{L^2}
    \ls M\ve_0\sum_{c\le a}\cE_1^c(t)+E_{|a|+1}(t).
\end{equation*}
Therefore,  for small $\ve_0>0$, we complete the proof of \eqref{WT:energy1}. \hfill $\square$

\section{On the analysis   of solutions away from the light cone}\label{sect4}

In this section, for the smooth solution $(v, f, g)$ of \eqref{axisymEuler}, we will
achieve the $L^\infty$ time-decay estimates away from the light cone for their space and time
derivatives through the structure of the hyperbolic system \eqref{axisymEuler}.

{\bf Lemma 4.1.} {\it For $a\in\N_0^2$ with $|a|\le N-3$, if
$(v, f, g)$ is  the solution of \eqref{axisymEuler}, then for small $\ve_0>0$, the following inequalities hold
\begin{eqnarray}
    |\chi_0(\p_r+\frac1r)\Gamma^af(t,r)|+|\chi_0\p_t\Gamma^af(t,r)|\ls
    <t>^{-1.99}\big[E_{|a|+2}(t)+\cX_{|a|+3}(t)+\cY_{|a|+2}(t)\big],\label{awayLC:f:pw}\\
    |\chi_0\p_r\Gamma^a\vv(t,r)|+|\chi_0\p_t\Gamma^a\vv(t,r)|\ls
    <t>^{-1.99}\big[E_{|a|+2}(t)+\cX_{|a|+3}(t)+\cY_{|a|+2}(t)\big],\label{awayLC:v:pw}\\
    |\p_t\Gamma^aG(t,r)|\ls
    <t>^{-1.99}\big[E_{|a|+2}(t)+\cX_{|a|+3}(t)+\cY_{|a|+2}(t)\big],\label{dtg:pw}\\
    |\p_t\Gamma^ag(t,r)|+|r\p_r\p_t\Gamma^aG(t,r)|\ls
    <t>^{-1.99}\big[E_{|a|+2}(t)+\cX_{|a|+3}(t)+\cY_{|a|+2}(t)\big].\label{dtg:pw'}
\end{eqnarray}
}

{\bf Proof.} At first, it is easy to find the following equations of $\vv, f$
\begin{equation}\label{axisymEuler0}
\left\{
\begin{aligned}
    &\p_t\vv-(\p_rf+\frac1rf)=\QQ_1,\\
    &\p_tf-\p_r\vv=\QQ_2,
\end{aligned}
\right.
\end{equation}
where
\begin{eqnarray}
  &&\QQ_1=: (\vv+G)(\p_rf+\frac1rf)-f\p_r\vv-f\p_rG-\p_tG,\label{QQ1}\\
  &&\QQ_2=: (\vv+G)\p_r\vv-f\p_rf.\label{QQ2}
\end{eqnarray}

According to \eqref{commutation},
we have the following  equations on the higher order derivatives of $\vv, f, g$
\begin{equation}\label{highordereqn0}
\left\{
\begin{aligned}
    &\p_t\Gamma^a\vv-(\p_r\Gamma^af+\frac1r\Gamma^af)=\tilde\Gamma^a\QQ_1,\\
    &\p_t\Gamma^af-\p_r\Gamma^a\vv=\tilde\Gamma^a\QQ_2,\\
    &\p_t\Gamma^ag=-\sum_{b+c=a}\sigma_{bc}\Gamma^bf(\p_r+\frac1r)\Gamma^cg.
\end{aligned}
\right.
\end{equation}

Similarly to \eqref{weighteqn1}, it is easy to get the following identities
\begin{equation}\label{weighteqn0}
\begin{array}{l}
\displaystyle (t^2-r^2)\p_r\Gamma^a\vv=tr\tilde\Gamma^a\QQ_1-t^2\tilde\Gamma^a\QQ_2
+tS\Gamma^af-rS\Gamma^a\vv+t\Gamma^af,\\
\displaystyle (t^2-r^2)\p_t\Gamma^a\vv=-r^2\tilde\Gamma^a\QQ_1+tr\tilde\Gamma^a\QQ_2
+tS\Gamma^a\vv-rS\Gamma^af-r\Gamma^af,\\
\displaystyle (t^2-r^2)\p_t\Gamma^af=tr\tilde\Gamma^a\QQ_1-r^2\tilde\Gamma^a\QQ_2
-rS\Gamma^a\vv+tS\Gamma^af+t\Gamma^af,\\
\displaystyle (t^2-r^2)(\p_r\Gamma^af+\frac1r\Gamma^af)=-t^2\tilde\Gamma^a\QQ_1
+tr\tilde\Gamma^a\QQ_2+tS\Gamma^a\vv-rS\Gamma^af-r\Gamma^af.
\end{array}
\end{equation}

Applying Leibniz's formula and \eqref{commutation},
we achieve
\begin{equation}\label{highorderQQ1}
    \tilde\Gamma^a\QQ_1=\sum_{b+c=a}\sigma_{bc}[\Gamma^b(\vv+G)(\p_r+\frac1r)\Gamma^cf
    -\Gamma^bf\p_r\Gamma^c\vv-\frac{\Gamma^bf}{r}r\p_r\Gamma^cG]-\p_t\Gamma^aG,
\end{equation}
and
\begin{equation}\label{highorderQQ2}
    \tilde\Gamma^a\QQ_2=\sum_{b+c=a}\sigma_{bc}[\Gamma^b(\vv+G)\p_r\Gamma^c\vv
    -\Gamma^bf(\p_r+\frac1r)\Gamma^cf+\frac1r\Gamma^bf\Gamma^cf].
\end{equation}

Next, we deal with $\p_t\Gamma^aG$ in \eqref{highorderQQ1}.
Taking the time-derivative $\p_t$ on two sides of \eqref{G:eqn1} yields
\begin{equation}\label{G:eqn2}
\begin{array}{l}
\displaystyle (1+v)\p_r\p_t\Gamma^aG+\sum_{b+c=a,c<a}
\sigma_{bc}(\Gamma^b\vv+\Gamma^bG)\p_r\p_t\Gamma^cG\\
\displaystyle =-\sum_{b+c=a}\sigma_{bc}(\p_t\Gamma^b\vv+\p_t\Gamma^bG)\p_r\Gamma^cG
-\frac2r\sum_{b+c=a}\frac{a!}{b!c!}\Gamma^bg\p_t\Gamma^cg.
\end{array}
\end{equation}
For the last term in \eqref{G:eqn2}, it is easy to get that
\begin{equation}\label{G:eqn3}
    \sum_{b+c=a}\frac{a!}{b!c!}\Gamma^bg\p_t\Gamma^cg=
    -\sum_{b+c+d=a}\frac{a!}{b!c!d!}\Gamma^bg\Gamma^cf(\p_r+\frac1r)\Gamma^dg.
\end{equation}
Therefore, we obtain the following inequality which is similar to \eqref{G:L1norm1}
\begin{equation}\label{G:L1norm3}
\begin{array}{l}
\displaystyle(1-\|\chi_0\vv\|_{L^\infty}-\|G\|_{L^\infty})\int_0^\infty|\p_r\p_t\Gamma^aG|dr\\
\displaystyle \ls\sum_{b+c=a,c<a}(\|\chi_0\Gamma^b\vv\|_{L^\infty}
+\|\Gamma^bG\|_{L^\infty})\int_0^\infty|\p_r\p_t\Gamma^cG|dr\\
\displaystyle \qquad+\sum_{b+c=a}(\|\chi_0\p_t\Gamma^b\vv\|_{L^\infty}
+\|\p_t\Gamma^bG\|_{L^\infty})\int_0^\infty|\p_r\Gamma^cG|dr\\
\displaystyle \qquad+\sum_{b+c+d=a}\|\frac1r\chi_0\Gamma^cf\|_{L^\infty}
\int_0^\infty|\Gamma^bg(\p_r+\frac1r)\Gamma^dg|dr.
\end{array}
\end{equation}
Here we point out that the second line in \eqref{G:L1norm3} does not appear if $|a|=0$.
Applying \eqref{hardyineq2} to the last term in \eqref{G:L1norm3}, we see that
\begin{equation}\label{g:hardy2}
\begin{array}{l}
\displaystyle \int_0^\infty|\Gamma^bg(\p_r+\frac1r)\Gamma^dg|dr\ls
\int_0^\infty\frac{|\Gamma^bg|^2}{r}dr+\int_0^\infty|(\p_r+\frac1r)\Gamma^dg|^2rdr\\
\displaystyle \ls\int_0^\infty|(\p_r+\frac1r)\Gamma^bg|^2rdr
+\int_0^\infty|(\p_r+\frac1r)\Gamma^dg|^2rdr.
\end{array}
\end{equation}

We now turn to treating the cross cutoff domain $\chi_0\chi_1\tilde\Gamma^a\QQ_1$ and $\chi_0\chi_1\tilde\Gamma^a\QQ_2$.
Direct computation yields
\begin{equation}\label{cross:cutoff}
\begin{split}
\displaystyle \chi_0\chi_1\tilde\Gamma^a\QQ_1=\sum_{b+c=a}\sigma_{bc}
[\chi_0\Gamma^b\vv\chi_1\p_r\Gamma^c(v+f)-\chi_1\Gamma^b(v+f)\chi_0\p_r\Gamma^c\vv
+\frac1r\chi_0\Gamma^b\vv\chi_1\Gamma^cf],\\
\displaystyle \chi_0\chi_1\tilde\Gamma^a\QQ_2=\sum_{b+c=a}\sigma_{bc}
[\chi_0\Gamma^b\vv\chi_1\p_r\Gamma^c(v+f)-\chi_1\Gamma^b(v+f)\chi_0\p_r\Gamma^cf].
\end{split}
\end{equation}

Collecting \eqref{weighteqn0}--\eqref{highorderQQ2}, \eqref{G:L1norm3}--\eqref{cross:cutoff}, and Lemma~2.5-2.6 with the method of induction,
we get \eqref{awayLC:f:pw}-\eqref{dtg:pw}.

Finally, we go to the proof of \eqref{dtg:pw'}.
For this purpose, at first, we need to treat $\|\p_t\Gamma^ag\|_{L^\infty}$.
It follows from the third equation in \eqref{highordereqn0} that
\begin{equation*}
    \|\p_t\Gamma^ag\|_{L^\infty}\ls\sum_{b+c=a}\|\Gamma^bf(\p_r+\frac1r)\Gamma^cg\|_{L^\infty}.
\end{equation*}
In addition, note that
\begin{equation}\label{dr:regular}
    r\p_r\Gamma^c=S\Gamma^c-t\p_t\Gamma^c.
\end{equation}
Then, by Lemma~2.5 and \eqref{awayLC:f:pw} we conclude
\begin{equation*}
\begin{array}{l}
\displaystyle\|\p_t\Gamma^ag\|_{L^\infty}\ls\sum_{b+c=a}\|\chi_0\frac1r\Gamma^bf\|_{L^\infty}
\big[\|(S+1)\Gamma^cg\|_{L^\infty}+\|t\p_t\Gamma^cg\|_{L^\infty}\big]\\
\displaystyle \ls <t>^{-1.99}\sum_{b+c=a}\big[E_{|b|+2}(t)+\cX_{|b|+3}(t)
+\cY_{|b|+2}(t)\big]\big[E_{|c|+2}(t)+\|t\p_t\Gamma^cg\|_{L^\infty}\big].
\end{array}
\end{equation*}
Due to the smallness of $\ve_0>0$, we deduce that
\begin{equation}\label{dtg:pw''}
\begin{array}{l}
\displaystyle \|\p_t\Gamma^ag\|_{L^\infty}\ls <t>^{-1.99}\sum_{b+c=a}E_{|c|+2}(t)
\big[E_{|b|+2}(t)+\cX_{|b|+3}(t)+\cY_{|b|+2}(t)\big]\\
\displaystyle \ls <t>^{-1.99}\big[E_{|a|+2}(t)+\cX_{|a|+3}(t)+\cY_{|a|+2}(t)\big].
\end{array}
\end{equation}

Finally, we deal with $\|r\p_r\p_t\Gamma^aG\|_{L^\infty}$.
From \eqref{G:eqn2}, one knows that
\begin{equation*}
\begin{array}{l}
\displaystyle \|r\p_r\p_t\Gamma^aG\|_{L^\infty}\ls\sum_{b+c=a,c<a}\|r\p_r\p_t\Gamma^cG\|_{L^\infty}
\big[\|\chi_0\Gamma^b\vv\|_{L^\infty}+\|\Gamma^bG\|_{L^\infty}\big]\\
\displaystyle \qquad+\sum_{b+c=a}\|r\p_r\Gamma^cG\|_{L^\infty}
\big[\|\chi_0\p_t\Gamma^b\vv\|_{L^\infty}+\|\p_t\Gamma^bG\|_{L^\infty}\big]
+\sum_{b+c=a}\|\Gamma^bg\p_t\Gamma^cg\|_{L^\infty}.
\end{array}
\end{equation*}
This, together with \eqref{awayLC:v:pw}-\eqref{dtg:pw}, \eqref{dtg:pw''} and Lemma~2.5-2.6, yields
\begin{equation*}
    \|r\p_r\p_t\Gamma^aG\|_{L^\infty}\ls <t>^{-1.99}\big[E_{|a|+2}(t)+\cX_{|a|+3}(t)+\cY_{|a|+2}(t)\big].
\end{equation*}
Therefore, we finish the proof of \eqref{dtg:pw'}. \hfill $\square$

\vskip 0.1 true cm

The following results will play an important role in the proof of Lemma~4.3 below.
\vskip 0.1 true cm

{\bf Lemma 4.2.} {\it
For $a,a'\in\N_0^2$ with $|a|\le N-3$ and $|a'|\le N-2$, if
$(v, f, g)$ is  the  solution of \eqref{axisymEuler}, then for small $\ve_0>0$,
the following inequalities hold

(1) $|\frac1r\Gamma^ag(t,r)|\ls\|(\p_r+\frac1r)\Gamma^ag(t,r)\|_{L^\infty}
    \ls E_{|a|+1}(t)+\cX_{|a|+2}(t)+\cY_{|a|+1}(t)+W_{|a|}(t);$

(2) $\|\p_r\Gamma^{a'}G\|_{L^\infty}\ls E_{|a'|+1}(t)+\cX_{|a'|+2}(t)+\cY_{|a'|+1}(t);$

(3) $\|\p_r\p_t\Gamma^aG\|_{L^\infty}\ls<t>^{-1.99}[E_{|a|+2}(t)+\cX_{|a|+3}(t)
    +\cY_{|a|+2}(t)+W_{|a|}(t)].$
}

\vskip 0.1 true cm

{\bf Proof.} First, we treat (1).
Similarly to (2) in Lemma~2.5, we achieve
\begin{equation*}
    \left|\frac1r\Gamma^ag(t,r)\right|=\frac1{r^2}\left|\int_0^r\p_r(r'\Gamma^ag(t,r'))dr'\right|
    \ls\|(\p_r+\frac1r)\Gamma^ag(t,r)\|_{L^\infty}.
\end{equation*}

Let $\phi(t,x)=\Gamma^aw(t,r)$ in \eqref{sobolev:W13}, then we obtain
\begin{equation*}
    \|\Gamma^aw(t,r)\|_{L^\infty}\ls\|\Gamma^aw(t,r)\|_{L^3}
    +\|\p_r\Gamma^aw(t,r)\|_{L^3}\ls W_{|a|}(t).
\end{equation*}

Recalling the definition of the specific vorticity $w$, we find that
\begin{equation*}
    \|(\p_r+\frac1r)\Gamma^ag(t,r)\|_{L^\infty}=\|\tilde\Gamma^a(\p_r+\frac1r)g(t,r)\|_{L^\infty}
    =\|\tilde\Gamma^a\Big(\frac{w}{1+\vv+G}\Big)\|_{L^\infty}.
\end{equation*}
Subsequently, we achieve (1).\\

Next, we focus on the estimate of $\|\p_r\Gamma^{a'}G\|_{L^\infty}$.
By making use of \eqref{G:eqn1}, we arrive at
\begin{equation*}
    \|\p_r\Gamma^{a'}G\|_{L^\infty}\ls\sum_{b+c=a',b<a'}\|\p_r\Gamma^cG\|_{L^\infty}
    \big[E_{|b|+1}(t)+\cX_{|b|+2}(t)+\cY_{|b|+1}(t)\big]+\sum_{b+c=a'}\|\frac1r\Gamma^bg\Gamma^cg\|_{L^\infty},
\end{equation*}
where
\begin{equation*}
    \sum_{b+c=a'}\|\frac1r\Gamma^bg\Gamma^cg\|_{L^\infty}\ls E_{|a'|+1}(t)
    \big[E_{[\frac{|a'|}{2}]+1}(t)+\cX_{[\frac{|a'|}{2}]+2}(t)+\cY_{[\frac{|a'|}{2}]+1}(t)
    +W_{[\frac{|a'|}{2}]}(t)\big]\ls E_{|a'|+1}(t).
\end{equation*}
Consequently, (2) is derived.\\

Similarly, \eqref{G:eqn2} implies
\begin{eqnarray*}
    \|\p_r\p_t\Gamma^aG\|_{L^\infty}\ls\sum_{b+c=a,c<a}\|\p_r\p_t\Gamma^cG\|_{L^\infty}
    \big[E_{|b|+1}(t)+\cX_{|b|+2}(t))+\cY_{|b|+1}(t)\big]\\
    +\sum_{b+c=a}\|\p_t\Gamma^b(\vv+G)\p_r\Gamma^cG\|_{L^\infty}+\sum_{b+c=a}\|\frac1r\Gamma^bg\p_t\Gamma^cg\|_{L^\infty},
\end{eqnarray*}
where
\begin{eqnarray*}
    &&\sum_{b+c=a}\|\frac1r\Gamma^bg\p_t\Gamma^cg\|_{L^\infty}\\
    &&\ls<t>^{-1.99}\sum_{b+c=a}\big[E_{|b|+1}(t)+\cX_{|b|+2}(t)+\cY_{|b|+1}(t)+W_{|b|}(t)\big]
    \big[E_{|c|+2}(t)+\cX_{|c|+3}(t)+\cY_{|c|+2}(t)\big].
\end{eqnarray*}
Note that $b+c=a$ and $|a|\le N-3$, either $|b|\le N-4$ or $|c|\le N-5$ holds.
Subsequently, we infer that
\begin{equation*}
    \sum_{b+c=a}\|\frac1r\Gamma^bg\p_t\Gamma^cg\|_{L^\infty}
    \ls<t>^{-1.99}[E_{|a|+2}(t)+\cX_{|a|+3}(t)+\cY_{|a|+2}(t)+W_{|a|}(t)].
\end{equation*}
Therefore, we obtain (3). \hfill $\square$

\vskip 0.1 true cm
{\bf Lemma 4.3.}  {\it
For any integer $k\in\N$ with $k\le N$ and $a\in\N_0^2$ with $|a|\le N-1$, if
$(v, f, g)$ is  the solution of \eqref{axisymEuler}, then  for small $\ve_0>0$,
the following weighted energy inequalities hold
\begin{eqnarray}
    &\cX_k(t)\ls E_k(t),\label{WT:energy0}\\
    &\|\p_t\Gamma^ag\|_{L^2}+\|r\p_r\p_t\Gamma^aG\|_{L^2}+\|\p_t\Gamma^aG\|_{L^2}\ls<t>^{-1}E_{|a|+1}(t).\label{dtg:L2norm}
\end{eqnarray}
}
\vskip 0.1 true cm
{\bf Proof.} We divide the proof of Lemma 4.3 into two parts.
\vskip 0.1 true cm
{\textbf{Part I: first order space derivatives estimates.}}

\vskip 0.1 true cm
For $a\in\N_0^2$ with $|a|\le N-1$, set
\begin{align*}
    \cE_0^a(t)&=:\|\chi_0\p_t\Gamma^a\vv(t,r)\|_{L^2}+\|\chi_0\p_r\Gamma^a\vv(t,r)\|_{L^2}\\
    & +\|\chi_0\p_t\Gamma^af(t,r)\|_{L^2}+\|\chi_0(\p_r+\frac1r)\Gamma^af(t,r)\|_{L^2}.
\end{align*}

According to \eqref{weighteqn0}, we have
\begin{equation}\label{WE01}
    \cE_0^a(t)\ls<t>^{-1}E_{|a|+1}(t)+\|\chi_0\tilde\Gamma^a\QQ_1\|_{L^2}
    +\|\chi_0\tilde\Gamma^a\QQ_2\|_{L^2}.
\end{equation}

At first, we deal with the term $\Gamma^b(\vv+G)(\p_r+\frac1r)\Gamma^cf$ in $\tilde\Gamma^a\QQ_1$.
If $|b|\le N-4$, we find that
\begin{eqnarray*}
    \|\chi_0\Gamma^b(\vv+G)(\p_r+\frac1r)\Gamma^cf\|_{L^2}\ls\cE_0^c(t)
    \big[\|\Gamma^bG\|_{L^\infty}+\|\chi_0\Gamma^b\vv\|_{L^\infty}+\|\chi_1\Gamma^bv\|_{L^\infty}\big]\\
    \ls\cE_0^c(t)\big[E_{|b|+1}(t)+\cX_{|b|+2}(t)+\cY_{|b|+2}(t)\big]\ls M\ve_0\cE_0^c(t),
\end{eqnarray*}
where we have used \eqref{assumption}, \eqref{v:pointwise} and Lemma~2.5-2.6.
For $|b|\ge N-3$, due to $|b|+|c|=|a|\le N-1$ and $N\ge N_0$ ,
one has $|c|\le2\le N-3$. Therefore, we achieve
\begin{equation}\label{WE01'}
\begin{array}{l}
    \|\chi_0\Gamma^b(\vv+G)(\p_r+\frac1r)\Gamma^cf\|_{L^2}\ls
    \|\Gamma^b(\vv+G)\|_{L^2}\|\chi_0(\p_r+\frac1r)\Gamma^cf\|_{L^\infty}\\
    \ls<t>^{-1.99}E_{|b|}(t)\big[E_{|c|+2}(t)+\cX_{|c|+3}(t)+\cY_{|c|+2}(t)\big]\ls<t>^{-1}E_{|b|}(t).
\end{array}
\end{equation}

The treatments for $\Gamma^bf\p_r\Gamma^c\vv, \Gamma^b(\vv+G)\p_r\Gamma^c\vv$
and $\Gamma^bf(\p_r+\frac1r)\Gamma^cf$ in $\tilde\Gamma^a\QQ_1$ and $\tilde\Gamma^a\QQ_2$ are analogous.
Then one has
\begin{equation}\label{WE02}
\begin{array}{l}
\displaystyle \|\chi_0\tilde\Gamma^a\QQ_1\|_{L^2}+\|\chi_0\tilde\Gamma^a\QQ_2\|_{L^2}
\ls M\ve_0\sum_{c\le a}\cE_0^c(t)+<t>^{-1}E_{|a|}(t)\\
\displaystyle \qquad+\sum_{b+c=a}\big[\|\chi_0\frac{\Gamma^bf}{r}r\p_r\Gamma^cG\|_{L^2}
+\|\chi_0\frac1r\Gamma^bf\Gamma^cf\|_{L^2}+\|\p_t\Gamma^aG\|_{L^2}\big].
\end{array}
\end{equation}
Now we deal with the terms in the last line of \eqref{WE02}.
By an analogous analysis for the multi-index $b$ and $c$ as in \eqref{WE01'}, we achieve
\begin{equation*}
\begin{array}{l}
\displaystyle \sum_{b+c=a}\big[\|\chi_0\frac{\Gamma^bf}{r}r\p_r\Gamma^cG\|_{L^2}
+\|\chi_0\frac1r\Gamma^bf\Gamma^cf\|_{L^2}\big]\\
\displaystyle \ls\sum_{b+c=a,|b|\le N-3}\|\chi_0\frac{\Gamma^bf}{r}\|_{L^\infty}
\big[\|r\p_r\Gamma^cG\|_{L^2}+\|\Gamma^cf\|_{L^2}\big]\\
\displaystyle \qquad+\sum_{b+c=a,|c|\le N-4}\|\chi_0\frac{\Gamma^bf}{r}\|_{L^2}
\big[\|r\p_r\Gamma^cG\|_{L^\infty}+\|\Gamma^cf\|_{L^\infty}\big].
\end{array}
\end{equation*}
Applying \eqref{hardyineq3} and Lemma~2.5--Lemma~4.1 yields
\begin{equation}\label{WE03}
\begin{array}{l}
\displaystyle \sum_{b+c=a}\big[\|\chi_0\frac{\Gamma^bf}{r}r\p_r\Gamma^cG\|_{L^2}
+\|\chi_0\frac1r\Gamma^bf\Gamma^cf\|_{L^2}\big]\\
\displaystyle \ls<t>^{-1.99}\sum_{b+c=a,|b|\le N-3}E_{|c|}(t)
\big[E_{|b|+2}(t)+\cX_{|b|+3}(t)+\cY_{|b|+2}(t)]\\
\displaystyle \qquad+\sum_{b+c=a,|c|\le N-4}
\big[\cE_0^b(t)+<t>^{-1}E_{|b|}(t)\big]\big[E_{|c|+1}(t)+\cX_{|c|+2}(t)+\cY_{|c|+1}(t)\big]\\
\displaystyle \ls M\ve_0\sum_{b\le a}\cE_0^b(t)+<t>^{-1}E_{|a|}(t).
\end{array}
\end{equation}

Next we turn our attention to the last term $\|\p_t\Gamma^aG\|_{L^2}$ in \eqref{WE02}.
Direct computation yields
\begin{equation*}
    \|\p_t\Gamma^aG\|^2_{L^2}=\int_0^\infty|\p_t\Gamma^aG|^2d\frac{r^2}{2}=
    \left|\int\p_t\Gamma^aG\p_r\p_t\Gamma^aGr^2dr\right|\le
    \|\p_t\Gamma^aG\|_{L^2}\|r\p_r\p_t\Gamma^aG\|_{L^2}.
\end{equation*}
Then we achieve
\begin{equation}\label{dtG:L2norm1}
    \|\p_t\Gamma^aG\|_{L^2}\le\|r\p_r\p_t\Gamma^aG\|_{L^2}.
\end{equation}

To deal with $\|r\p_r\p_t\Gamma^aG\|_{L^2}$ in \eqref{dtG:L2norm1},
taking $L^2$ norms on the both side of \eqref{G:eqn2} derives
\begin{equation}\label{dtG:L2norm2}
\begin{array}{l}
\displaystyle \|r\p_r\p_t\Gamma^aG\|_{L^2}\ls\sum_{b+c=a,c<a}
\|(\Gamma^b\vv+\Gamma^bG)r\p_r\p_t\Gamma^cG\|_{L^2}\\
\displaystyle\qquad+\sum_{b+c=a}\big[\|(\p_t\Gamma^b\vv+\p_t\Gamma^bG)r\p_r\Gamma^cG\|_{L^2}
+\|\Gamma^bg\p_t\Gamma^cg\|_{L^2}\big].
\end{array}
\end{equation}
Note that $\chi_1\equiv0$ on the support of $G$.
Then similar to \eqref{WE02} and \eqref{WE03}, we obtain
\begin{equation}\label{dtG:L2norm3}
\begin{array}{l}
\displaystyle \sum_{b+c=a,c<a}\|(\Gamma^b\vv+\Gamma^bG)r\p_r\p_t\Gamma^cG\|_{L^2}
\ls\sum_{b+c=a,|c|\le N-3}\big[\|\Gamma^b\vv\|_{L^2}+\|\Gamma^bG\|_{L^2}\big]
\|r\p_r\p_t\Gamma^cG\|_{L^\infty}\\
\displaystyle \qquad+\sum_{b+c=a,c<a,|b|\le N-4}\|r\p_r\p_t\Gamma^cG\|_{L^2}
\big[\|\chi_0\Gamma^b\vv\|_{L^\infty}+\|\Gamma^bG\|_{L^\infty}\big]\\
\displaystyle \ls M\ve_0\sum_{c\le a}\|r\p_r\p_t\Gamma^cG\|_{L^2}+<t>^{-1}E_{|a|}(t),
\end{array}
\end{equation}
and
\begin{equation}\label{dtG:L2norm4}
\begin{array}{l}
\displaystyle \sum_{b+c=a}\|(\p_t\Gamma^b\vv+\p_t\Gamma^bG)r\p_r\Gamma^cG\|_{L^2}\\
\displaystyle \ls\sum_{b+c=a,|c|\le N-4}\|r\p_r\Gamma^cG\|_{L^\infty}
\big[\|\chi_0\p_t\Gamma^b\vv\|_{L^2}+\|\p_t\Gamma^bG\|_{L^2}\big]\\
\displaystyle \qquad+\sum_{b+c=a,|b|\le N-3}\|r\p_r\Gamma^cG\|_{L^2}
\big[\|\chi_0\p_t\Gamma^b\vv\|_{L^\infty}+\|\p_t\Gamma^bG\|_{L^\infty}\big]\\
\displaystyle \ls M\ve_0\sum_{b\le a}\big[\cE_0^b(t)+\|r\p_r\p_t\Gamma^bG\|_{L^2}\big]
+<t>^{-1}E_{|a|}(t).
\end{array}
\end{equation}
In addition, according to the third equation in \eqref{highordereqn0} and \eqref{dr:regular},
we deduce that
\begin{equation*}
\begin{array}{l}
\displaystyle \|\p_t\Gamma^cg\|_{L^2}\ls\sum_{d+e=c,|d|\le N-3}
\|\chi_0\frac1r\Gamma^df\|_{L^\infty}
\big[\|(S+1)\Gamma^eg\|_{L^2}+\|t\p_t\Gamma^eg\|_{L^2}\big]\\
\displaystyle \quad+\sum_{d+e=c,|e|\le N-3}\|\chi_0\frac1r\Gamma^df\|_{L^2}
\big[\|(S+1)\Gamma^eg\|_{L^\infty}+\|t\p_t\Gamma^eg\|_{L^\infty}\big]\\
\displaystyle \ls<t>^{-1.99}\sum_{d+e=c,|d|\le N-3}
\big[E_{|d|+2}(t)+\cX_{|d|+3}(t)+\cY_{|d|+2}(t)\big]\big[E_{|e|+1}(t)+<t>\|\p_t\Gamma^eg\|_{L^2}\big]\\
\displaystyle \quad+\sum_{\substack{d+e=c,\\|e|\le N-3}}
\big[\cE_0^d(t)+<t>^{-1}E_{|d|}(t)\big]
\big\{E_{|e|+1}(t)+<t>^{-0.99}[E_{|e|+2}(t)+\cX_{|e|+3}(t)+\cY_{|e|+2}(t)]\big\},
\end{array}
\end{equation*}
where we have used \eqref{hardyineq3}, Lemma~2.5 and Lemma~4.1. Consequently,
\begin{equation*}\label{dtg:L2norm1}
    \|\p_t\Gamma^cg\|_{L^2}\ls M\ve_0\cE_0^c(t)+<t>^{-1}E_{|c|+1}(t).
\end{equation*}
This yields
\begin{equation}\label{dtg:L2norm2}
    \sum_{b+c=a}\|\Gamma^bg\p_t\Gamma^cg\|_{L^2}\ls
    M\ve_0\sum_{c\le a}\cE_0^c(t)+<t>^{-1}E_{|a|+1}(t).
\end{equation}
Inserting \eqref{dtG:L2norm3}, \eqref{dtG:L2norm4} and \eqref{dtg:L2norm2} into \eqref{dtG:L2norm2}, we arrive at
\begin{equation*}
    \|\p_t\Gamma^aG\|_{L^2}\ls\|r\p_r\p_t\Gamma^aG\|_{L^2}
    \ls M\ve_0\sum_{b\le a}\cE_0^b(t)+<t>^{-1}E_{|a|+1}(t).
\end{equation*}
This, together with \eqref{WE01}, \eqref{WE02} and \eqref{WE03}, yields that for all $|a|\le N-1$,
\begin{equation}\label{WE04}
    \cE_0^a(t)+\|\p_t\Gamma^ag\|_{L^2}+\|r\p_r\p_t\Gamma^aG\|_{L^2}+\|\p_t\Gamma^aG\|_{L^2}\ls<t>^{-1}E_{|a|+1}(t).
\end{equation}

\vskip 0.1 true cm

{\textbf{Part II: second order space derivatives estimates.}}
\vskip 0.1 true cm

Let $|a'|\le N-2$.
From \eqref{highordereqn0}, we easily get
\begin{equation*}
\begin{array}{l}
\displaystyle   \p_r(\p_r+\frac1r)\Gamma^{a'}f=\p_r\p_t\Gamma^{a'}\vv-\p_r\tilde\Gamma^{a'}\QQ_1\\
\displaystyle   (\p_r+\frac1r)\p_r\Gamma^{a'}\vv=(\p_r+\frac1r)\p_t\Gamma^{a'}f-(\p_r+\frac1r)\tilde\Gamma^{a'}\QQ_2,
\end{array}
\end{equation*}
which yields
\begin{equation}\label{WE05}
\begin{array}{l}
\displaystyle \cV^{a'}(t)\defeq\|\chi_0\p_r(\p_r+\frac1r)\Gamma^{a'}f\|_{L^2}
+\|\chi_0(\p_r+\frac1r)\p_r\Gamma^{a'}v\|_{L^2}\\
\displaystyle\ls\|\chi_0\p_r\p_t\Gamma^{a'}\vv\|_{L^2}+\|\chi_0(\p_r+\frac1r)\p_t\Gamma^{a'}f\|_{L^2}
+\|\chi_0\p_r\tilde\Gamma^{a'}\QQ_1\|_{L^2}+\|\chi_0(\p_r+\frac1r)\tilde\Gamma^{a'}\QQ_2\|_{L^2}.
\end{array}
\end{equation}

We now treat the last two nonlinear terms in \eqref{WE05}.
It is not hard to check that
\begin{equation}\label{drQQ1}
\begin{array}{l}
\displaystyle \p_r\tilde\Gamma^{a'}\QQ_1=\sum_{b+c=a'}\sigma_{bc}
[\Gamma^b(\vv+G)\p_r(\p_r+\frac1r)\Gamma^cf-\Gamma^bf(\p_r+\frac1r)\p_r\Gamma^c\vv\\
\displaystyle \qquad+2\Gamma^bf\frac1r\p_r\Gamma^c\vv+\frac1r\Gamma^bf\p_r\Gamma^cG
-\Gamma^bf\p_r\p_r\Gamma^cG]-\p_r\p_t\Gamma^{a'}G,
\end{array}
\end{equation}
and
\begin{equation}\label{drQQ2}
\begin{array}{l}
\displaystyle (\p_r+\frac1r)\tilde\Gamma^{a'}\QQ_2=\sum_{b+c=a'}
\sigma_{bc}[\Gamma^b(\vv+G)(\p_r+\frac1r)\p_r\Gamma^c\vv-\Gamma^bf\p_r(\p_r+\frac1r)\Gamma^cf\\
\displaystyle \qquad+r\p_r\Gamma^bG\frac1r\p_r\Gamma^c\vv+\p_r\Gamma^b\vv\p_r\Gamma^c\vv
-(\p_r+\frac1r)\Gamma^bf(\p_r+\frac1r)\Gamma^cf+2\frac{\Gamma^bf\p_r\Gamma^cf}{r}].
\end{array}
\end{equation}
Applying \eqref{hardyineq3} to the first term $\frac1r\p_r\Gamma^c\vv$
in the second line of \eqref{drQQ1}, we achieve
\begin{equation}\label{drQQ1:L2}
\begin{array}{l}
\displaystyle \|\chi_0\p_r\tilde\Gamma^{a'}\QQ_1\|_{L^2}\ls\sum_{b+c=a'}\big[\cV^c(t)
+<t>^{-1}E_{|c|+1}(t)\big]\big[E_{|b|}(t)+\cX_{|b|+2}(t)+\cY_{|b|+1}(t)\big]\\
\displaystyle \qquad+\sum_{b+c=a'}\big[\|\chi_0\frac1r\Gamma^bf\p_r\Gamma^cG\|_{L^2}
+\|\chi_0\Gamma^bf\p_r\p_r\Gamma^cG\|_{L^2}\big]+\|\p_r\p_t\Gamma^{a'}G\|_{L^2}.
\end{array}
\end{equation}
Next we start to treat each term in the last line of \eqref{drQQ1:L2}.

\vskip 0.2 true cm

\textbf{A. Estimate of $\|\chi_0\frac1r\Gamma^bf\p_r\Gamma^cG\|_{L^2}$.}

\vskip 0.2 true cm

For $|b|\le N-3$, from \eqref{G:eqn1}, we see that
\begin{equation*}
\begin{array}{l}
\displaystyle \|\chi_0\frac1r\Gamma^bf\p_r\Gamma^cG\|_{L^2}\ls\sum_{d+e=c,e<c}
\|\chi_0\frac1r\Gamma^bf\Gamma^d(\vv+G)\p_r\Gamma^eG\|_{L^2}
+\sum_{d+e=c}\|\chi_0\frac1r\Gamma^bf\Gamma^dg\frac1r\Gamma^eg\|_{L^2}\\
\displaystyle \ls\sum_{d+e=c,e<c}E_{|d|}(t)\|\chi_0\frac1r\Gamma^bf\p_r\Gamma^eG\|_{L^\infty}
+\sum_{d+e=c}\|\chi_0\frac1r\Gamma^bf\Gamma^dg\|_{L^\infty}\|\frac1r\Gamma^eg\|_{L^2},
\end{array}
\end{equation*}
where
\begin{eqnarray*}
    &&\sum_{b+d+e=a'}\|\chi_0\frac1r\Gamma^bf\Gamma^dg\|_{L^\infty}\|\frac1r\Gamma^eg\|_{L^2}\\
    &&\ls<t>^{-1.99}\sum_{b+c=a',|b|\le N-3}E_{[\frac{|c|}{2}]+1}E_{|c|+1}
    \big[E_{|b|+2}(t)+\cX_{|b|+3}(t)+\cY_{|b|+2}(t)\big]\\
    &&\ls<t>^{-1}E_{|a'|+1}(t).
\end{eqnarray*}
This yields
\begin{equation*}
    \sum_{b+c=a',|b|\le N-3}\|\chi_0\frac1r\Gamma^bf\p_r\Gamma^cG\|_{L^2}
    \ls<t>^{-1}E_{|a'|+1}(t).
\end{equation*}

Next, we deal with the case of $|b|\ge N-2$.
Applying \eqref{hardyineq3} to $\|\chi_0\frac1r\Gamma^bf\|_{L^2}$
and the first inequality in Lemma~4.2 to $\|\frac1r\Gamma^eg\|_{L^\infty}$, we obtain
\begin{eqnarray*}
    &&\sum_{d+e=c}\|\chi_0\frac1r\Gamma^bf\|_{L^2}\|\Gamma^dg\|_{L^\infty}\|\frac1r\Gamma^eg\|_{L^\infty}\\
    &&\ls\sum_{d+e=c}E_{|d|+1}(t)\big[\cE_0^b(t)+<t>^{-1}E_{|b|}(t)\big]
    \big[E_{|e|+1}(t)+\cX_{|e|+2}(t)+\cY_{|e|+1}(t)+W_{|e|}(t)\big]\\
    &&\ls M\ve_0\cE_0^b(t)+<t>^{-1}E_{|b|}(t).
\end{eqnarray*}
Therefore,
\begin{equation}\label{drQQ1:L2:1}
    \sum_{b+c=a'}\|\chi_0\frac1r\Gamma^bf\p_r\Gamma^cG\|_{L^2}
    \ls M\ve_0\sum_{b\le a}\cE_0^b(t)+<t>^{-1}E_{|a'|+1}(t).
\end{equation}

\vskip 0.1 true cm
\textbf{B. Estimate of $\|\chi_0\Gamma^bf\p_r\p_r\Gamma^cG\|_{L^2}$.}
\vskip 0.1 true cm

Taking the space derivative $\p_r$ on two sides of \eqref{G:eqn1}, we find that
\begin{eqnarray*}
    &&(1+v)\p_r\p_r\Gamma^cG+\sum_{d+e=c,e<c}\sigma_{de}\Gamma^d(\vv+G)\p_r\p_r\Gamma^eG\\
    &&=-\sum_{d+e=c}\sigma_{de}\big[\p_r\Gamma^d(\vv+G)\p_r\Gamma^eG
    +\frac2r\Gamma^dg\p_r\Gamma^eg-\frac{1}{r^2}\Gamma^dg\Gamma^eg\big].
\end{eqnarray*}

For $|b|\le N-3$, we arrive at
\begin{equation*}
\begin{array}{l}
\displaystyle \|\chi_0\Gamma^bf\p_r\p_r\Gamma^cG\|_{L^2}\ls\sum_{d+e=c,e<c}
\|\chi_0\Gamma^d(\vv+G)\Gamma^bf\p_r\p_r\Gamma^eG\|_{L^2}\\
\displaystyle \qquad+\|\chi_0\frac1r\Gamma^bf\|_{L^\infty}\sum_{d+e=c}\big[\|\chi_0\p_r\Gamma^d(\vv+G)r\p_r\Gamma^eG\|_{L^2}
+\|\Gamma^dg(\p_r+\frac1r)\Gamma^eg\|_{L^2}+\|\Gamma^dg\frac1r\Gamma^eg\|_{L^2}\big],
\end{array}
\end{equation*}
where
\begin{equation*}
    \sum_{d+e=c,e<c}\|\chi_0\Gamma^d(\vv+G)\Gamma^bf\p_r\p_r\Gamma^eG\|_{L^2}\\
    \ls\sum_{d+e=c,|e|\le N-3}E_{|d|}(t)\|\frac1r\Gamma^bfr\p_r\p_r\Gamma^eG\|_{L^\infty}.
\end{equation*}
With the help of \eqref{dr:regular}, we deduce that
\begin{equation*}
    r\p_r\p_r\Gamma^eG=\p_rS\Gamma^eG-\p_r\Gamma^eG-t\p_r\p_t\Gamma^eG,
\end{equation*}
which yields
\begin{equation*}
    \|r\p_r\p_r\Gamma^eG\|_{L^\infty}\ls E_{|e|+2}(t)+\cX_{|e|+3}(t)+\cY_{|e|+2}(t)+W_{|e|}(t).
\end{equation*}
Consequently, we arrive at
\begin{equation*}
    \|\chi_0\Gamma^bf\p_r\p_r\Gamma^cG\|_{L^2}\ls<t>^{-1}E_{|c|+1}(t).
\end{equation*}

While for $|b|\ge N-2$, we have
\begin{equation*}
\begin{array}{l}
\displaystyle \|\chi_0\Gamma^bf\p_r\p_r\Gamma^cG\|_{L^2}\ls\sum_{d+e=c,e<c}
\|\chi_0\Gamma^bf\p_r\p_r\Gamma^eG\|_{L^2}\big[E_{|d|+1}(t)+\cX_{|d|+2}(t)\big]\\
\displaystyle \qquad+\|\chi_0\frac1r\Gamma^bf\|_{L^2}\sum_{d+e=c}
\big[\|\chi_0\p_r\Gamma^d(\vv+G)r\p_r\Gamma^eG\|_{L^\infty}
+\|\Gamma^dg(\p_r+\frac1r)\Gamma^eg\|_{L^\infty}
+\|\frac1r\Gamma^dg\Gamma^eg\|_{L^\infty}\big]\\
\displaystyle \ls M\ve_0\sum_{e<c}\|\chi_0\Gamma^bf\p_r\p_r\Gamma^eG\|_{L^2}
+M\ve_0\cE_0^b(t)+<t>^{-1}E_{|b|}(t).
\end{array}
\end{equation*}
Thus,
\begin{equation}\label{drQQ1:L2:2}
    \sum_{b+c=a'}\|\chi_0\Gamma^bf\p_r\p_r\Gamma^cG\|_{L^2}
    \ls M\ve_0\sum_{b\le a}\cE_0^b(t)+<t>^{-1}E_{|a'|+1}(t).
\end{equation}

\vskip 0.1 true cm
\textbf{C. Estimate of $\|\p_r\p_t\Gamma^{a'}G\|_{L^2}$.}
\vskip 0.1 true cm

Applying similar analysis to \eqref{G:eqn2}, we arrive at
\begin{equation*}
\begin{array}{l}
\displaystyle \|\p_r\p_t\Gamma^{a'}G\|_{L^2}\ls\sum_{b+c=a',c<a'}
\|\chi_0\Gamma^b(\vv+G)\p_r\p_t\Gamma^cG\|_{L^2}\\
\displaystyle \qquad+\sum_{b+c=a',|b|\le N-3}
\big[\|\frac1r\Gamma^bg\|_{L^\infty}\|\p_t\Gamma^cg\|_{L^2}
+\|\chi_0\p_t\Gamma^b(\vv+G)\|_{L^\infty}\|\p_r\Gamma^cG\|_{L^2}\big]\\
\displaystyle \qquad+\sum_{b+c=a',|b|\ge N-2}
\big[\|\frac1r\Gamma^bg\|_{L^2}\|\p_t\Gamma^cg\|_{L^\infty}
+\|\chi_0\p_t\Gamma^b(\vv+G)\|_{L^2}\|\p_r\Gamma^cG\|_{L^\infty}\big]\\
\displaystyle \ls M\ve_0\sum_{b\le a'}\big[\cE_0^b(t)+\|\p_r\p_t\Gamma^bG\|_{L^2}\big]
+<t>^{-1}E_{|a'|+1}(t).
\end{array}
\end{equation*}
This, together with \eqref{drQQ1:L2}, \eqref{drQQ1:L2:1} and \eqref{drQQ1:L2:2}, yields
\begin{equation}\label{drQQ1:L2'}
\begin{array}{l}
\displaystyle \|\chi_0\p_r\tilde\Gamma^{a'}\QQ_1\|_{L^2}\ls\sum_{b+c=a'}
\big[\cV^c(t)+<t>^{-1}E_{|c|+1}(t)\big]\big[E_{|b|}(t)+\cX_{|b|+2}(t)+\cY_{|b|+1}(t)\big]\\
\displaystyle \qquad+M\ve_0\sum_{b\le a}\cE_0^b(t)+<t>^{-1}E_{|a'|+1}(t).
\end{array}
\end{equation}
With an analogous analysis to $\|\chi_0\p_r\tilde\Gamma^{a'}\QQ_1\|_{L^2}$, we achieve
\begin{equation}\label{drQQ2:L2}
\begin{array}{l}
\displaystyle \|\chi^2_0(\p_r+\frac1r)\tilde\Gamma^{a'}\QQ_2\|_{L^2}\ls\sum_{b+c=a'}
\big[\cV^c(t)+<t>^{-1}E_{|c|+1}(t)\big]\big[E_{|b|}(t)+\cX_{|b|+2}(t)+\cY_{|b|+1}(t)\big]\\
\displaystyle \qquad+M\ve_0\sum_{b\le a}\cE_0^b(t)+<t>^{-1}E_{|a'|+1}(t).
\end{array}
\end{equation}

In addition, due to
\begin{eqnarray*}
    \chi_0\chi_1(\p_r+\frac1r)\tilde\Gamma^{a'}\QQ_2=\sum_{b+c=a'}\sigma_{bc}
    [\chi_1\Gamma^bv\chi_0(\p_r+\frac1r)\p_r\Gamma^c\vv-\chi_1\Gamma^bf\chi_0\p_r(\p_r+\frac1r)\Gamma^cf\\
    +\chi_1\p_r\Gamma^b(v+f)\chi_0\p_r\Gamma^c(\vv-f)-\frac{1}{r^2}\chi_0\chi_1\Gamma^bf\Gamma^cf],
\end{eqnarray*}
then by Lemma~2.2 and Lemma~3.1, we arrive at
\begin{equation}\label{drQQ2:L2'}
    \|\chi_0\chi_1(\p_r+\frac1r)\tilde\Gamma^{a'}\QQ_2\|_{L^2}\ls<t>^{-\frac12}
    \sum_{b+c=a'}\big[\cV^c(t)+<t>^{-1}E_{|c|+1}(t)\big]\big[E_{|b|}(t)+\cY_{|b|+2}(t)\big].
\end{equation}
Collecting \eqref{WE05} and \eqref{drQQ1:L2'}-\eqref{drQQ2:L2'} yields
\begin{equation*}
\begin{array}{l}
\displaystyle \cV^{a'}(t)\ls\sum_{b+c=a'}
\big[\cV^c(t)+<t>^{-1}E_{|c|+1}(t)\big]\big[E_{|b|}(t)+\cX_{|b|+2}(t)+\cY_{|b|+1}(t)\big]\\
\displaystyle \qquad+M\ve_0\sum_{b\le a}\cE_0^b(t)+<t>^{-1}E_{|a'|+1}(t).
\end{array}
\end{equation*}
Then it leads to \eqref{WT:energy0} by the method of induction on $|a'|$.

Thus, by collecting Part I and Part II,  Lemma 4.3 is proved. \hfill $\square$

\section{Elementary energy estimates  of solutions}\label{sect5}

With the help of \eqref{WT:energy1} and \eqref{WT:energy0}, Lemma~2.2,~3.1,~3.2 and Lemma~2.5,~2.6,~4.1 can be summarized as follows.

\vskip 0.1 true cm

{\bf Lemma 5.1.} {\it For $a,a',b\in\N_0^2$ with $|a|\le N-1$, $|a'|\le N-2$ and $|b|\le N-3$,
if $(v, f, g)$ is  the solution of \eqref{axisymEuler}, then for small $\ve_0>0$,
the following inequalities hold
\begin{equation*}
\begin{array}{l}
\displaystyle |<r-t>^\frac12\chi_1\Gamma^af(t,r)|+|<r-t>^\frac12\chi_1\Gamma^av(t,r)|
\ls<t>^{-\frac12}E_{|a|+1}(t),\\
\displaystyle |\chi_1\p_r\Gamma^{a'}(v+f)|\ls<t>^{-\frac32}E_{|a'|+2}(t),
\quad\|\chi_1\p_r\Gamma^a(v+f)\|_{L^2}\ls<t>^{-1}E_{|a|+1}(t),\\
\displaystyle |<r-t>^\frac32\chi_1\p_r\Gamma^bf(t,r)|
+|<r-t>^\frac32\chi_1\p_t\Gamma^bf(t,r)|\ls<t>^{-\frac12}E_{|b|+2}(t),\\
\displaystyle |<r-t>^\frac32\chi_1\p_r\Gamma^bv(t,r)|
+|<r-t>^\frac32\chi_1\p_t\Gamma^bv(t,r)|\ls<t>^{-\frac12}E_{|b|+2}(t).
\end{array}
\end{equation*}
}

{\bf Lemma 5.2.} {\it For $a,a',b\in\N_0^2$ with $|a|\le N-1$, $|a'|\le N-2$ and $|b|\le N-3$,
if $(v, f, g)$ is  the solution of \eqref{axisymEuler}, then for small $\ve_0>0$,
the following inequalities hold
\begin{equation*}
\begin{array}{l}
\displaystyle |\Gamma^ag(t,r)|\ls E_{|a|+1}(t),
\quad|\Gamma^{a'}G(t,r)|+|r\p_r\Gamma^{a'}G(t,r)| \ls E_{|a'|+2}(t),\\
\displaystyle |\chi_0\Gamma^af(t,r)|\ls<t>^{-1}E_{|a|+1}(t),
\quad|\chi_0\Gamma^{a'}\vv(t,r)|\ls<t>^{-0.99}E_{|a'|+2}(t),\\
\displaystyle |\chi_0(\p_r+\frac1r)\Gamma^bf(t,r)|+|\chi_0\p_t\Gamma^bf(t,r)|
+|\frac1r\chi_0\Gamma^bf(t,r)|\ls<t>^{-1.99}E_{|b|+3}(t),\\
\displaystyle |\chi_0\p_r\Gamma^b\vv(t,r)|+|\chi_0\p_t\Gamma^b\vv(t,r)|
\ls<t>^{-1.99}E_{|b|+3}(t),\\
\displaystyle |\p_t\Gamma^bg(t,r)|+|r\p_r\p_t\Gamma^bG(t,r)|+|\p_t\Gamma^bG(t,r)|\ls<t>^{-1.99}E_{|b|+3}(t).
\end{array}
\end{equation*}
}

Based on Lemma 5.1 and Lemma 5.2, we establish the following results.

\vskip 0.1 true cm

{\bf Lemma 5.3.} {\it For $a\in\N_0^2$, if
$(v, f, g)$ is  the solution of \eqref{axisymEuler}, then  for small $\ve_0>0$,
the following elementary energy inequality holds
\begin{equation}\label{energyineq}
\begin{array}{l}
\displaystyle \|\Gamma^av(t,r)\|^2_{L^2}+\|\Gamma^af(t,r)\|^2_{L^2}
+\int_0^t\|<r-t'>^{-\frac58}\Gamma^a(v+f)\|^2_{L^2}dt'\\
\displaystyle \ls\|\Gamma^av(0,r)\|^2_{L^2}+\|\Gamma^af(0,r)\|^2_{L^2}
+\int_0^t<t'>^{-1.99}E_3(t')E^2_{|a|}(t')dt'\\
\displaystyle \qquad+\left|\int_0^t\int(I-e^qr\p_r\vv\Gamma^av\Gamma^af)drdt'\right|,
\end{array}
\end{equation}
where
\begin{equation}\label{I:def}
\begin{array}{l}
\displaystyle I=:\frac{1+\vv}{1+v}\Big\{-q'e^qrv\Gamma^av\Gamma^af
-\frac12e^qrf\p_r(|\Gamma^av|^2+|\Gamma^af|^2)+e^qr|\Gamma^av|^2(\p_r+\frac1r)f\\
\displaystyle \qquad-e^qr|\Gamma^af|^2\p_rf+e^q\Gamma^af\Gamma^a(g^2)
+\sum_{b+c=a,b<a,c<a}\sigma_{bc}e^q\Gamma^av\Gamma^bv\Gamma^cf\\
\displaystyle \qquad+\sum_{\substack{b+c=a,\\b<a,c<a}}\sigma_{bc}e^qr\Big[\p_r\Gamma^cv
(\Gamma^af\Gamma^bv-\Gamma^av\Gamma^bf)+\p_r\Gamma^cf(\Gamma^av\Gamma^bv-\Gamma^af\Gamma^bf)\Big]\Big\}
\end{array}
\end{equation}
and the smooth function $q=q(r-t)$ satisfies $q'(s)=<s>^{-5/4}$
and $q(\infty)=0$.
}
\vskip 0.1 true cm
{\bf Remark 5.1.} {\it The multiplier function $e^{q(r-t)}$ in \eqref{energyineq} is called the ``ghost weight"
by S.~Alinhac in \cite{Alinhac01}.}
\vskip 0.1 true cm

{\bf Proof.} Multiplying \eqref{highordereqn} by $e^qr\Gamma^av$ and $e^qr\Gamma^af$ respectively,
we achieve
\begin{equation}\label{energyidentity}
\begin{array}{l}
\displaystyle \frac12\p_t[e^qr(|\Gamma^av|^2+|\Gamma^af|^2)]
-\p_r(e^qr\Gamma^av\Gamma^af)+\frac12q'e^qr|\Gamma^a(v+f)|^2\\
\displaystyle =e^qr(\Gamma^av\tilde\Gamma^aQ_1+\Gamma^af\tilde\Gamma^aQ_2).
\end{array}
\end{equation}

We now treat the terms containing top-order derivative $\Gamma^av$ and $\Gamma^af$ in the right hand side
of \eqref{energyidentity}. To this end, we rewrite \eqref{highorderQ1}-\eqref{highorderQ2} as follows:
\begin{equation*}
    \tilde\Gamma^aQ_1=v(\p_r+\frac1r)\Gamma^af-f\p_r\Gamma^av+\Gamma^av(\p_r+\frac1r)f
    -\Gamma^af\p_rv+\sum_{\substack{b+c=a,\\b<a,c<a}}\sigma_{bc}[\Gamma^bv(\p_r+\frac1r)\Gamma^cf-\Gamma^bf\p_r\Gamma^cv]
\end{equation*}
and
\begin{equation*}
    \tilde\Gamma^aQ_2=v\p_r\Gamma^av-f\p_r\Gamma^af+\Gamma^av\p_rv-\Gamma^af\p_rf
    +\frac1r\Gamma^a(g^2)+\sum_{\substack{b+c=a,\\b<a,c<a}}\sigma_{bc}[\Gamma^bv\p_r\Gamma^cv-\Gamma^bf\p_r\Gamma^cf].
\end{equation*}

It follows from direct computation that
\begin{eqnarray*}
    &&e^qr(\Gamma^av\tilde\Gamma^aQ_1+\Gamma^af\tilde\Gamma^aQ_2)\\
    &&=v\p_r(e^qr\Gamma^av\Gamma^af)-q'e^qrv\Gamma^av\Gamma^af
    -\frac12e^qrf\p_r(|\Gamma^av|^2+|\Gamma^af|^2)+e^qr|\Gamma^av|^2(\p_r+\frac1r)f\\
    &&\quad-e^qr|\Gamma^af|^2\p_rf+e^q\Gamma^af\Gamma^a(g^2)
    +\sum_{b+c=a,b<a,c<a}\sigma_{bc}e^q\Gamma^av\Gamma^bv\Gamma^cf\\
    &&\quad+\sum_{b+c=a,b<a,c<a}\sigma_{bc}e^qr\Big[\p_r\Gamma^cv(\Gamma^af\Gamma^bv-\Gamma^av\Gamma^bf)
    +\p_r\Gamma^cf(\Gamma^av\Gamma^bv-\Gamma^af\Gamma^bf)\Big].
\end{eqnarray*}
Inserting this into \eqref{energyidentity} yields
\begin{equation*}
    \frac12\p_t[e^qr(|\Gamma^av|^2+|\Gamma^af|^2)]+\frac12q'e^qr|\Gamma^a(v+f)|^2
    =(1+v)\p_r(e^qr\Gamma^av\Gamma^af)+\frac{1+v}{1+\vv}I.
\end{equation*}
Multiplying this identity by $\frac{1+\vv}{1+v}$
and integrating it over $[0,t]\times[0,+\infty)$, we have
\begin{equation}\label{energyineq1}
\begin{split}
&\|\Gamma^av(t,r)\|^2_{L^2}+\|\Gamma^af(t,r)\|^2_{L^2}
+\int_0^t\|<r-t'>^{-\frac58}\Gamma^a(v+f)\|^2_{L^2}dt'\\
&\ls\left|\int_0^t\int(I-e^qr\p_r\vv\Gamma^av\Gamma^af)drdt'\right|
+\int_0^t\int\Big|\p_t\Big(\frac{1+\vv}{1+v}\Big)\Big|[|\Gamma^av|^2+|\Gamma^af|^2]rdrdt'.
\end{split}
\end{equation}
Note that $G\equiv0$ on $\supp\chi_1$, thus from Lemma~5.2 we know
\begin{equation}\label{energyineq2}
\begin{split}
\Big|\p_t\Big(\frac{1+\vv}{1+v}\Big)\Big|\ls\Big|\chi_0\p_t\Big(\frac{1+\vv}{1+v}\Big)\Big|
\ls|\p_tG|+|\chi_0\p_t\vv|\ls<t'>^{-1.99}E_3(t').
\end{split}
\end{equation}
Substituting \eqref{energyineq2} into \eqref{energyineq1} derives \eqref{energyineq}.
\hfill $\square$

\vskip 0.1 true cm
In the subsequent $\S 6$ and $\S 7$, we will treat the term 
\begin{equation}\label{H0}
\begin{split}
\int_0^t\int(I-e^qr\p_r\vv\Gamma^av\Gamma^af)drdt'
\end{split}
\end{equation}
in the right hand side of \eqref{energyineq} near the light cone
and away from the light cone, respectively.

\section{The treatment  of \eqref{H0} near the light cone}\label{sect6}

Based on the results in Section 5, we now treat \eqref{H0} near the light cone
for the smooth solution $(v, f, g)$ of \eqref{axisymEuler}.
\vskip 0.1 true cm

{\bf Lemma 6.1.} {\it For $a\in\N_0^2$ with $|a|\le N$, if $(v, f, g)$ is  the
solution of \eqref{axisymEuler}, then   for small $\ve_0>0$, the following energy inequality near the light cone holds
\begin{equation}\label{nearLC:ineq}
\begin{array}{l}
\displaystyle \left|\int_0^t\int\chi_1(I-e^qr\p_r\vv\Gamma^av\Gamma^af)drdt'\right|\\
\displaystyle \ls\int_0^t<t'>^{-\frac32}E_N(t')E^2_{|a|}(t')dt'\\
+\left\{
\begin{aligned}
    &\sum_{b\le a}\int_0^tE_{N-2}(t')\|\sqrt{q'}\Gamma^b(v+f)\|_{L^2}dt'
    +\int_0^t<t'>^{-1}E_{N-2}(t')E^2_{|a|}(t')dt',\;|a|\le N,\\
    &\int_0^t<t'>^{-\frac32}E_{N}(t')E_{|a|}^2(t')dt',\;|a|\le N-2,
\end{aligned}
\right.
\end{array}
\end{equation}
where the smooth function $q=q(r-t)$ in \eqref{nearLC:ineq} satisfies $q'(s)=<s>^{-5/4}$
and $q(\infty)=0$.
}
\vskip 0.1 true cm
{\bf Proof.} Because of $g\equiv0$ and $G\equiv0$ on $\supp\chi_1$, we find that $v\equiv\vv$.
Recalling the definition of $I$ in \eqref{I:def}, we see that
\begin{eqnarray*}
    &&\chi_1I=\chi_1\Big\{-q'e^qrv\Gamma^av\Gamma^af-\frac12e^qrf\p_r(|\Gamma^av|^2+|\Gamma^af|^2)
    +e^qr\p_rf(|\Gamma^av|^2-|\Gamma^af|^2)\\
    &&\qquad\quad+e^qf|\Gamma^av|^2\Gamma^af+\sum_{b+c=a,b<a,c<a}\sigma_{bc}e^qr\Gamma^av\Gamma^bv\Gamma^cf\\
    &&\qquad\quad+\sum_{\substack{b+c=a,\\b<a,c<a}}\sigma_{bc}e^qr[\p_r\Gamma^c(v+f)
    (\Gamma^af\Gamma^bv-\Gamma^av\Gamma^bf)+\Gamma^a(v-f)\Gamma^b(v+f)\p_r\Gamma^cf]\Big\}.
\end{eqnarray*}
Set
\begin{equation*}
    \chi_1(I-e^qr\p_rv\Gamma^av\Gamma^af)=\chi_1(I_{11}+I_{12}+I_{13}+I_{14}),
\end{equation*}
where
\begin{eqnarray*}
   &&I_{11}=: -e^qr\p_rv\Gamma^av\Gamma^af-q'e^qrv\Gamma^av\Gamma^af
   -\frac12e^qrf\p_r(|\Gamma^av|^2+|\Gamma^af|^2)\\
   &&\qquad\quad+e^qr\p_rf(|\Gamma^av|^2-|\Gamma^af|^2),\\
   &&I_{12}=: e^qf|\Gamma^av|^2+\sum_{b+c=a,b<a,c<a}\sigma_{bc}e^q
   \Gamma^av\Gamma^bv\Gamma^cf,\\
   &&I_{13}=: \sum_{b+c=a,b<a,c<a}\sigma_{bc}e^qr\Gamma^a(v-f)\Gamma^b(v+f)\p_r\Gamma^cf\\
   &&I_{14}=: \sum_{b+c=a,b<a,c<a}\sigma_{bc}e^qr\p_r\Gamma^c(v+f)
   (\Gamma^af\Gamma^bv-\Gamma^av\Gamma^bf).
\end{eqnarray*}
Applying integration by parts to the third term in $I_{11}$, then the resulting inequality can be divided into three parts as follows:
\begin{equation}\label{I11ineq}
    \left|\int_0^t\int\chi_1I_{11}drdt'\right|\ls I_{11}^1+I_{11}^2+I_{11}^3,
\end{equation}
where
\begin{eqnarray*}
   &&I_{11}^1=:\int_0^t\int(<t'>^{-1}|\chi'_1|+\chi_1\frac1r)|f|(|\Gamma^av|^2+|\Gamma^af|^2)rdrdt',\\
   &&I_{11}^2=:\int_0^t\int\chi_1\big[q'|(v+f)|+|\p_r(v+f)|\big]|\Gamma^av\Gamma^af|rdrdt',\\
   &&I_{11}^3=:\int_0^t\int\chi_1\big[q'|f|+|\p_rf|\big]|\Gamma^a(v+f)|(|\Gamma^av|+|\Gamma^af|)rdrdt'.
\end{eqnarray*}

In view of Lemma~5.1, we achieve
\begin{equation}\label{I111+12ineq}
    \left|\int_0^t\int\chi_1I_{12}drdt'\right|+I_{11}^1
    \ls\int_0^t<t'>^{-\frac32}E_{[\frac{|a|}{2}]+1}(t')E^2_{|a|}(t')dt'.
\end{equation}

By using \eqref{goodLinfty'} for the term $q'|(v+f)|$ in $I_{11}^2$, we derive
\begin{equation}\label{I112ineq}
    I_{11}^2\ls\int_0^t<t'>^{-\frac32}E_2(t')E^2_{|a|}(t')dt'.
\end{equation}

For the higher-order energy with $N-1\le|a|\le N$ in $I_{11}^3$ and $|c|\le N-4$ in $I_{13}$,
applying the Cauchy-Schwartz inequality and Lemma~5.1 to $I_{11}^3$ and $I_{13}$ implies
\begin{equation}\label{I113ineq:higher}
\begin{array}{l}
\displaystyle I_{11}^3+\sum_{b<a,|c|\le N-4}\left|\int_0^t\int
\chi_1e^q\Gamma^a(v-f)\Gamma^b(v+f)\p_r\Gamma^cfrdrdt'\right|\\
\displaystyle \ls\sum_{b\le a,|c|\le N-4}\int_0^tE_{|c|+2}(t')
\|\sqrt{q'}\Gamma^b(v+f)\|_{L^2}dt'+\int_0^t<t'>^{-1}E_{|c|+2}(t')E^2_{|a|}(t')dt'\\
\displaystyle \ls\sum_{b\le a}\int_0^tE_{N-2}(t')\|\sqrt{q'}\Gamma^b(v+f)\|_{L^2}dt'
+\int_0^t<t'>^{-1}E_{N-2}(t')E^2_{|a|}(t')dt'.
\end{array}
\end{equation}
For the lower-order energy with $|a|\le N-2$ in $I_{11}^3$, by using \eqref{goodLinfty'} again
for the term $\Gamma^a(v+f)$, we then have
\begin{equation*}
    I_{11}^3\ls\int_0^t\|<r-t'>^{-1}\chi_1\Gamma^a(v+f)\|_{L^\infty}\big[\|f\|_{L^2}
    +\|<r-t'>\p_rf\|_{L^2}\big]E_{|a|}(t'),
\end{equation*}
where
\begin{equation*}
    \|<r-t'>\p_rf\|_{L^2}\ls\|<r-t'>\chi_1\p_rf\|_{L^2}
    +\|<r-t'>\chi_0\p_rf\|_{L^2}\ls E_2(t').
\end{equation*}
This treatment can also be applied to $I_{13}$ for $|c|\ge N-3$. Therefore,we achieve
\begin{equation}\label{I113ineq:lower}
    I_{11}^3+\sum_{\substack{b+c=a,\\b<a,|c|\ge N-3}}\left|\int_0^t\int
    \chi_1e^q\Gamma^a(v-f)\Gamma^b(v+f)\p_r\Gamma^cfrdrdt'\right|
    \ls\int_0^t<t'>^{-\frac32}E_N(t')E_{|a|}^2(t')dt'.
\end{equation}

Finally, we turn our attention to the last term $I_{14}$.
Notice that $\chi_1\Gamma^bv=\chi_1\Gamma^b\vv$ holds. Then it is not hard to find that
\begin{eqnarray*}
    \left|\int_0^t\int\chi_1I_{14}drdt'\right|\ls\int_0^tE_{|a|}(t')\big[\sum_{|c|<[\frac{|a|}{2}]}
    \|\chi_1\p_r\Gamma^c(v+f)\|_{L^\infty}(\|\Gamma^bv\|_{L^2}+\|\Gamma^bf\|_{L^2})\\
    +\sum_{c<a,|b|\le[\frac{|a|}{2}]}\|\chi_1\p_r\Gamma^c(v+f)\|_{L^2}
    (\|\Gamma^b\vv\|_{L^\infty}+\|\Gamma^bf\|_{L^\infty})\big],
\end{eqnarray*}
where
\begin{equation*}
\begin{array}{l}
    \|\Gamma^b\vv\|_{L^\infty}\ls\|\chi_1\Gamma^bv\|_{L^\infty}
    +\|\chi_0\Gamma^b\vv\|_{L^\infty}\le<t'>^{-\frac12}E_{|b|+2}(t').
\end{array}
\end{equation*}
Thus,
\begin{equation}\label{I14ineq}
\begin{array}{l}
\displaystyle \left|\int_0^t\int\chi_1I_{14}drdt'\right|\\
\displaystyle \ls\int_0^t<t'>^{-\frac32}E_{|a|}(t')\big[\sum_{|c|<[\frac{|a|}{2}]}
E_{|b|}(t')E_{|c|+2}(t')+\sum_{c<a,|b|\le[\frac{|a|}{2}]}E_{|b|+2}(t')E_{|c|+1}(t')\big]\\
\displaystyle \ls\int_0^t<t'>^{-\frac32}E_{[\frac{|a|}{2}]+2}(t')E^2_{|a|}(t')dt'.
\end{array}
\end{equation}

Collecting \eqref{I11ineq}--\eqref{I14ineq} yields \eqref{nearLC:ineq}. \hfill $\square$

\section{The treatment  of \eqref{H0} away from the light cone and estimate of $w$}\label{sect7}

In this section, at first we treat \eqref{H0} away from the light cone
for the smooth solution $(v, f, g)$ of \eqref{axisymEuler}, and subsequently take the
estimate on $w$.

\vskip 0.1 true cm
{\bf Lemma 7.1.} {\it For  $a\in\N_0^2$ with $|a|\le N$,  if
$(v, f, g)$ is  the  solution of \eqref{axisymEuler}, then  for small $\ve_0>0$,
the following energy inequality away from the light cone holds
\begin{equation}\label{awayLC:ineq}
    \left|\int_0^t\int\chi_0(I-e^qr\p_r\vv\Gamma^av\Gamma^af)drdt'\right|
    \ls E_{\max\{3,[\frac{|a|}{2}]+2\}}(t)E^2_{|a|}(t)+\cQ_a(t),
\end{equation}
where
\begin{equation*}
\begin{array}{l}
\displaystyle \cQ_a(t)=:\int_0^t<t'>^{-\frac54}E^2_{|a|}(t')[E_N(t')+W_{N-3}(t')]dt'\\
\displaystyle \qquad\quad +\left\{\begin{aligned}
    &\int_0^t<t'>^{-1}E_{N-2}(t')E^2_{|a|}(t')dt',\quad|a|\le N,\\
    &\int_0^t<t'>^{-1.99}E_N(t')E^2_{|a|}(t')dt',\quad
    |a|\le N-2,
\end{aligned}
\right.
\end{array}
\end{equation*}
and the smooth function $q=q(r-t)$ in \eqref{awayLC:ineq} satisfies $q'(s)=<s>^{-5/4}$
and $q(\infty)=0$.
}

\vskip 0.1 true cm
{\bf Proof.} At first, it is easy to find that
\begin{equation}\label{I00ineq}
    \left|\int_0^t\int\chi_0e^qr\p_r\vv\Gamma^av\Gamma^afdrdt'\right|
    \ls\int_0^t<t'>^{-1.99}E_3(t')E^2_{|a|}(t')dt'.
\end{equation}

Next we focus on the treatment of $\chi_0I$ in \eqref{awayLC:ineq}. Let
\begin{equation}\label{I0:decomposition}
    \frac{1+v}{1+\vv}I=I_{01}+I_{02}+I_{03}+I_{04},
\end{equation}
where
\begin{eqnarray*}
   &&I_{01}=: -q'e^qrv\Gamma^av\Gamma^af+e^qr|\Gamma^av|^2(\p_rf+\frac1rf)
   -e^qr|\Gamma^af|^2\p_rf,\\
   &&I_{02}=: -\frac12e^qrf\p_r(|\Gamma^av|^2+|\Gamma^af|^2),\\
   &&I_{03}=: \sum_{b+c=a,b<a,c<a}\sigma_{bc}e^qr[(\p_r+\frac1r)\Gamma^cf
   (\Gamma^av\Gamma^bv-\Gamma^af\Gamma^bf)\\
   &&\qquad\quad+\p_r\Gamma^c\vv(\Gamma^af\Gamma^bv-\Gamma^av\Gamma^bf)
   +\frac{\Gamma^bf}{r}(\Gamma^af\Gamma^cf-r\p_r\Gamma^cG\Gamma^av)],\\
   &&I_{04}=: e^q\Gamma^af\Gamma^a(g^2)
   +\sum_{b+c=a,b<a,c<a}\sigma_{bc}e^qr\Gamma^af\Gamma^bv\p_r\Gamma^cG.
\end{eqnarray*}

By using Lemma~5.2 to $I_{01}$ directly, we see that
\begin{equation}\label{I01ineq}
    \left|\int_0^t\int\chi_0\frac{1+\vv}{1+v}I_{01}drdt'\right|
    \ls\int_0^t<t'>^{-\frac54}E_3(t')E^2_{|a|}(t')dt'\ls\cQ_a(t).
\end{equation}

Applying the integration by parts to $I_{02}$ derives
\begin{equation*}
    \left|\int_0^t\int\chi_0\frac{1+\vv}{1+v}I_{02}drdt'\right|\ls
    \left|\int_0^t\int\p_r(\chi_0e^qrf\frac{1+\vv}{1+v})[|\Gamma^av|^2+|\Gamma^af|^2]drdt'\right|,
\end{equation*}
where
\begin{eqnarray*}
    \frac1r|\p_r(\chi_0 e^qrf\frac{1+\vv}{1+v})|\ls|<t'>^{-1}\chi'_0f|
    +|q'\chi_0f|+|\chi_0(\p_rf+\frac1rf)|&\\
    +|\chi_0\frac{f}{r}r\p_rG|+|\chi_0f\p_r\vv|\ls<t'>^{-\frac32}E_3(t').
\end{eqnarray*}
Consequently, we achieve
\begin{equation}\label{I02ineq}
    \left|\int_0^t\int\chi_0\frac{1+\vv}{1+v}I_{02}drdt'\right|
    \ls\int_0^t<t'>^{-\frac32}E_3(t')E^2_{|a|}(t')dt'\ls\cQ_a(t).
\end{equation}

For the terms in $I_{03}$, it suffices only to deal with the term $\Gamma^af\Gamma^bv\p_r\Gamma^c\vv$
since the other terms can be analogously treated.
For the higher-energy with $N-1\le|a|\le N$, we arrive at
\begin{equation*}
\begin{array}{l}
\displaystyle \sum_{b+c=a,b<a,c<a}\|\chi_0\Gamma^bv\p_r\Gamma^c\vv\|_{L^2}\\
\displaystyle \ls\sum_{c<a,|b|\le N-4}\|\Gamma^bv\|_{L^\infty}
\|\chi_0\p_r\Gamma^c\vv\|_{L^2}
+\sum_{|c|\le N-3}\|\Gamma^bv\|_{L^2}\|\chi_0\p_r\Gamma^c\vv\|_{L^\infty}\\
\displaystyle \ls<t'>^{-1}\sum_{c<a,|b|\le N-4}E_{|b|+2}(t')E_{|c|+1}(t')
+<t'>^{-1.99}\sum_{|c|\le N-3}E_{|b|}(t')E_{|c|+3}(t')\\
\displaystyle \ls<t'>^{-1}E_{|a|}(t')E_{N-2}(t')
+<t'>^{-1.99}E_{|a|}(t')E_N(t').
\end{array}
\end{equation*}
While for the lower-energy with $|a|\le N-2$, we easily get
\begin{equation*}
\begin{array}{l}
\displaystyle \sum_{b+c=a,b<a,c<a}\|\chi_0\Gamma^bv\p_r\Gamma^c\vv\|_{L^2}
\displaystyle \ls\sum_{b+c=a,b<a,c<a}\|\Gamma^bv\|_{L^2}
\|\chi_0\p_r\Gamma^c\vv\|_{L^\infty}\\
\displaystyle \ls<t'>^{-1.99}\sum_{b+c=a,b<a,c<a}E_{|b|}(t')E_{|c|+3}(t')
\ls<t'>^{-1.99}E_{|a|}(t')E_N(t').
\end{array}
\end{equation*}
Therefore, we conclude
\begin{equation}\label{I03ineq}
\begin{array}{l}
\displaystyle \left|\int_0^t\int\chi_0\frac{1+\vv}{1+v}I_{03}drdt'\right|
\ls\int_0^t<t'>^{-1.99}E_N(t')E^2_{|a|}(t')dt'\\
\displaystyle \qquad+\left\{\begin{aligned}
    &\int_0^t<t'>^{-1}E^2_{|a|}(t')E_{N-2}(t')dt',\quad|a|\le N,\\
    &\int_0^t<t'>^{-1.99}E_{|a|}^2(t')E_N(t')dt',
    \quad|a|\le N-2,
\end{aligned}
\right.\\
\displaystyle \ls\cQ_a(t).
\end{array}
\end{equation}

Finally, we focus on the treatment of $I_{04}$.

If $\Gamma^a=\p_t^{a_1}S^{a_2}$ with $a_1\ge1$,
let $\p_t\Gamma^{a'}=\Gamma^a$ with $|a'|=|a|-1\le N-1$.
It is easy to find that
\begin{equation*}
    \Gamma^af\Gamma^a(g^2)=\p_t\Gamma^{a'}f\p_t\Gamma^{a'}(g^2) =\sum_{b+c=a'}\sigma_{bc}\p_t\Gamma^{a'}f\p_t(\Gamma^bg\Gamma^cg)
    =2\sum_{b+c=a'}\sigma_{bc}\p_t\Gamma^{a'}f\p_t\Gamma^bg\Gamma^cg.
\end{equation*}
Subsequently, we achieve
\begin{equation*}
\begin{array}{l}
\displaystyle\sum_{b+c=a'}\int\left|\chi_0\p_t\Gamma^{a'}f\p_t\Gamma^bg\Gamma^cg\right|dr\\
\displaystyle \ls\|\chi_0\p_t\Gamma^{a'}f\|_{L^2}
\big[\sum_{|b|\le N-3}\|\p_t\Gamma^bg\|_{L^\infty}\|\frac1r\Gamma^cg\|_{L^2}
+\sum_{|c|\le1}\|\p_t\Gamma^bg\|_{L^2}\|\frac1r\Gamma^cg\|_{L^\infty}\big]\\
\displaystyle \ls<t'>^{-1}E_{|a|}(t')
\big\{<t'>^{-1.99}\sum_{|b|\le N-3}E_{|b|+3}(t')E_{|c|+2}(t')\\
\displaystyle \qquad +<t'>^{-1}\sum_{|c|\le1}E_{|b|+1}(t')[E_{|c|+2}(t')+W_{|c|}(t')]\big\}\\
\displaystyle \ls<t'>^{-2}E^2_{|a|}(t')[E_N(t')+W_{N-3}(t')].
\end{array}
\end{equation*}
Analogously, for the other terms in $I_{04}$, one has
\begin{eqnarray*}
    &&\sum_{b+c=a,b<a,c<a}\sigma_{bc}\p_r\Gamma^cG\Gamma^af\Gamma^bv\\
    &&=\sum_{b+c=a',b<a'}\sigma_{bc}\p_r\Gamma^cG\p_t\Gamma^{a'}f\p_t\Gamma^bv
    +\sum_{b+c=a',c<a'}\sigma_{bc}\p_r\p_t\Gamma^cG\p_t\Gamma^{a'}f\Gamma^bv,
\end{eqnarray*}
which implies
\begin{equation*}
\begin{array}{l}
\displaystyle \sum_{b+c=a,b<a,c<a}
\int\left|\chi_0e^q\Gamma^af\Gamma^bv\p_r\Gamma^cG\right|rdr\\
\displaystyle \ls\|\chi_0\p_t\Gamma^{a'}f\|_{L^2}
\big[\sum_{|b|\le N-3}\|\chi_0\p_t\Gamma^bv\|_{L^\infty}\|\p_r\Gamma^cG\|_{L^2}
+\sum_{b<a,|c|\le2}\|\chi_0\p_t\Gamma^bv\|_{L^2}\|\p_r\Gamma^cG\|_{L^\infty}\\
\displaystyle \qquad\qquad+\sum_{c<a,|b|\le N-2}
\|\chi_0\Gamma^bv\|_{L^\infty}\|\p_r\p_t\Gamma^cG\|_{L^2}
+\sum_{|c|\le1}\|\Gamma^bv\|_{L^2}\|\p_r\p_t\Gamma^cG\|_{L^\infty}\big]\\
\displaystyle \ls<t'>^{-2}E_{|a|}(t')\big\{\sum_{|b|\le N-3}E_{|b|+3}(t')E_{|c|}(t')
+\sum_{b<a,|c|\le2}E_{|b|+1}(t')[E_{|c|+2}(t')+W_{|c|}(t')]\\
\displaystyle \qquad\qquad+\sum_{c<a,|b|\le N-2}E_{|b|+2}(t')E_{|c|+1}(t')
+\sum_{|c|\le1}E_{|b|}(t')[E_{|c|+3}(t')+W_{|c|}(t')]\big\}\\
\displaystyle \ls<t'>^{-2}E^2_{|a|}(t')[E_N(t')+W_{N-3}(t')],
\end{array}
\end{equation*}
where we have used Lemma~4.2. Consequently, we achieve that for $\Gamma^a=\p_t^{a_1}S^{a_2}$ with $a_1\ge1$
\begin{equation}\label{I04ineq:dt}
    \left|\int_0^t\int\chi_0\frac{1+\vv}{1+v}I_{04}drdt'\right|\ls\cQ_a(t).
\end{equation}

Next, we treat the case of $\Gamma^a=S^{a_2}=S^l$ and $l=|a|=a_2\ge0$.
When $l=|a|=0$, the second term $\ds\sum_{b+c=a,b<a,c<a}\sigma_{bc}e^qr\Gamma^af\Gamma^bv\p_r\Gamma^cG$ in $I_{04}$ does not appear.
While $l\ge1$, rewrite the scaling operator as $S^lf=(t'\p_t+r\p_r)S^{l-1}f$.
Then applying the integration by parts with respect to $t'$,
we arrive at that for all $m=1,\cdots,l-1$,
\begin{equation*}
\begin{array}{l}
\displaystyle \left|\int_0^t\int\chi_0e^q\frac{1+\vv}{1+v}S^lfS^mv\p_rS^{l-m}Grdrdt'\right|\\
\displaystyle \ls\int_0^t\int\left|\chi_0\p_rS^{l-1}fS^mvr\p_rS^{l-m}G\right|rdrdt'
+\int\left|\chi_0t\frac{S^{l-1}f}{r}S^mvr\p_rS^{l-m}G\right|rdr\\
\displaystyle \quad+\int_0^t\int\left|\p_t[\chi_0t'e^q\frac{1+\vv}{1+v}S^mvr\p_rS^{l-m}G]
\frac{S^{l-1}f}{r}\right|rdrdt'\\
\displaystyle \ls E_{[\frac{|a|}{2}]+2}(t)E^2_{|a|}(t)\\
\displaystyle \quad +\int_0^t<t'>^{-\frac54}E_N(t')E^2_{|a|}(t')dt'+\left\{\begin{aligned}
    &\int_0^t<t'>^{-1}E_{N-2}(t')E^2_{|a|}(t')dt',\quad|a|\le N,\\
    &\int_0^t<t'>^{-1.99}E_{|a|+2}(t')E^2_{|a|}(t')dt',\quad|a|\le N-2,
\end{aligned}
\right.\\
\displaystyle \ls E_{[\frac{|a|}{2}]+2}(t)E^2_{|a|}(t)+\cQ_a(t).
\end{array}
\end{equation*}

Next, we deal with the first term $S^lfS^l(g^2)$ in $I_{04}$.
When $l=|a|=0$, we achieve
\begin{equation}\label{I04ineq:a=0}
    \left|\int_0^t\int\chi_0\frac{1+\vv}{1+v}e^qfg^2drdt'\right|
    \ls\int_0^t\|\frac1r\chi_0f\|_{L^\infty}E^2_0(t')dt'\ls\int_0^tE_3(t')E^2_0(t')dt'.
\end{equation}
Similarly to the treatment for the second term in $I_{04}$, we have that
\begin{equation}\label{I04ineq:S}
\begin{array}{l}
\displaystyle \left|\int_0^t\int\chi_0\frac{1+\vv}{1+v}e^qS^lfS^l(g^2)drdt'\right|\\
\displaystyle \ls
\left|\int_0^t\int\chi_0t'\frac{1+\vv}{1+v}e^qg\p_tS^lgS^{l-1}fdrdt'\right|
+\int\left|\chi_0tS^l(g^2)\frac{S^{l-1}f}{r}\right|rdr\\
\displaystyle \qquad+\left|\int_0^t\int\p_t[\chi_0t'\frac{1+\vv}{1+v}e^q]S^l(g^2)\frac{S^{l-1}f}{r}rdrdt'\right|
+\int_0^t\int\left|\chi_0\p_rS^{l-1}fS^l(g^2)\right|rdrdt'\\
\displaystyle \qquad+\sum_{m=1}^{l-1}\int_0^t\int\left|\chi_0t'S^{l-m}g\p_tS^mg\frac{S^{l-1}f}{r}\right|rdrdt'
=:I_{04}^1+I_{04}^2+I_{04}^3+I_{04}^4+I_{04}^5.
\end{array}
\end{equation}
Note that $\p_tS^lg=-S^lf(\p_rg+\frac1rg)-f(\p_r+\frac1r)S^lg-\ds\sum_{m=1}^{l-1}S^{l-m}f(\p_r+\frac1r)S^mg$.
Substituting this into \eqref{I04ineq:S} yields
\begin{equation}\label{I041ineq}
\begin{array}{l}
\displaystyle I_{04}^1\ls
\left|\int_0^t\int\chi_0t'\frac{1+\vv}{1+v}e^qg(\p_rg+\frac1rg)S^lfS^{l-1}fdrdt'\right|\\
\displaystyle \qquad+\left|\int_0^t\int\chi_0t'\frac{1+\vv}{1+v}e^qg(\p_r+\frac1r)S^lgfS^{l-1}fdrdt'\right|\\ \displaystyle \qquad+\sum_{m=1}^{l-1}\int_0^t\int\left|\chi_0t'g(r\p_r+1)S^mg
\frac{S^{l-m}f}{r}\frac{S^{l-1}f}{r}\right|rdrdt'\\
\qquad =:I_{04}^{11}+I_{04}^{12}+I_{04}^{13}.
\end{array}
\end{equation}
It is easy to see that $S^lfS^{l-1}f=\frac12t'\p_t(|S^{l-1}f|^2)+r\p_rS^{l-1}fS^{l-1}f$.
Then we derive
\begin{equation*}
\begin{array}{l}
\displaystyle I_{04}^{11}\ls
\left|\int_0^t\int\p_t[\chi_0t'^2\frac{1+\vv}{1+v}e^qg(r\p_rg+g)]|\frac{S^{l-1}f}{r}|^2rdrdt'\right|\\
\displaystyle \qquad+\int\left|\chi_0t^2g(r\p_rg+g)|\frac{S^{l-1}f}{r}|^2\right|rdr
+\int_0^t\int\left|\chi_0t'g(r\p_rg+g)\p_rS^{l-1}f\frac{S^{l-1}f}{r}\right|rdrdt'.
\end{array}
\end{equation*}

Thanks to \eqref{dr:regular}, we deduce that
\begin{equation*}
    r\p_r\p_tg=-r\p_r[f\p_rg+\frac1rfg]=-\frac{f}{r}r\p_rSg
    +t\frac{f}{r}r\p_r\p_tg-\p_rf(r\p_rg+g)+\frac1rfg.
\end{equation*}
This, together with
\begin{equation*}
    \|r\p_rS^mg\|_{L^\infty}\ls\|S^{m+1}g\|_{L^\infty}+t\|\p_tS^mg\|_{L^\infty}\ls E_{m+2}(t')+<t'>^{-0.99}E_{m+3}(t'),
\end{equation*}
yields
\begin{equation*}
    \|r\p_r\p_tg\|_{L^\infty}\ls<t'>^{-1.99}E_3(t')E_4(t').
\end{equation*}
Therefore, we achieve
\begin{equation}\label{I0411ineq}
    I_{04}^{11}\ls E_3(t)E^2_{|a|}(t)+\cQ_a(t).
\end{equation}

For $I_{04}^{12}$, applying the integration by parts with respect to $r$, we arrive at
\begin{equation}\label{I0412ineq}
\begin{array}{l}
\displaystyle I_{04}^{12}\ls
\int_0^t\int\left|\chi_0t'\frac{1+\vv}{1+v}gS^lg\frac{\p_r(fS^{l-1}f)}{r}\right|rdrdt'\\
\displaystyle \qquad+\int_0^t\int\left|\chi_0t'gS^lg\frac{f}{r}\frac{S^{l-1}f}{r}\right|rdrdt'
+\left|\int_0^t\int r\p_r[\chi_0\frac{1+\vv}{1+v}e^qg]t'S^lg\frac{f}{r}\frac{S^{l-1}f}{r}rdrdt'\right|\\
\displaystyle \ls \int_0^t<t'>^{-1.99}E_3(t')E^2_{|a|}(t')dt'.
\end{array}
\end{equation}

For the other terms $I_{04}^{13}, I_{04}^{2}, I_{04}^{3}, I_{04}^{4}, I_{04}^{5}$,
applying the same analysis on the multi-index as in the above  we conclude
\begin{equation}\label{I0413ineq}
\begin{array}{l}
\displaystyle I_{04}^{13}+I_{04}^{2}+I_{04}^{3}+I_{04}^{4}+I_{04}^{5}
\displaystyle \ls E_{[\frac{|a|}{2}]+2}(t)E^2_{|a|}(t)+\cQ_a(t).
\end{array}
\end{equation}

Collecting \eqref{I00ineq}--\eqref{I0413ineq}, then \eqref{awayLC:ineq} is proved.
\hfill $\square$

\vskip 0.1 true cm
Next we estimate $g$ and $w$.

{\bf Lemma 7.2.} {\it For $a,a',b\in\N_0^2$ with $|a|\le N$, $|a'|\le N-1$ and $|b|\le N-3$,
if $(v, f, g)$ is  the solution of \eqref{axisymEuler}, then  for small $\ve_0>0$,
the following energy inequalities for $g$ and specific vorticity $w$ hold
\begin{eqnarray}
    \|\Gamma^ag(t,r)\|^2_{L^2}\ls\|\Gamma^ag(0,r)\|^2_{L^2}
    +E_{\max\{3,[\frac{|a|}{2}]+2\}}(t)E^2_{|a|}(t)+\cQ_a(t),\label{g:ineq}\\
    \|\Gamma^{a'}w(t,r)\|^2_{L^2}\ls\|\Gamma^{a'}w(0,r)\|^2_{L^2}
    +\tilde\cQ_{a'}(t),\label{curl:L2norm}\\
    \|\Gamma^bw(t,r)\|^3_{L^3}\ls\|\Gamma^bw(0,r)\|^3_{L^3}
    +\int_0^t<t'>^{-1.99}E_{|b|+3}(t')W^3_{|b|}(t')dt',\label{curl:L3norm}\\
    \|\p_r\Gamma^bw(t,r)\|^3_{L^3}\ls\|\p_r\Gamma^{b}w(0,r)\|^3_{L^3}
    +\int_0^t<t'>^{-1.99}E_{|b|+3}(t')W^3_{|b|}(t')dt',\label{curl:W13norm}
\end{eqnarray}
where
\begin{equation*}
\begin{array}{l}
\displaystyle \tilde\cQ_{a'}(t)=:
\left\{\begin{aligned}
    &\int_0^t<t'>^{-1}E_{N-2}(t')E^2_{|a'|+1}(t')dt',\quad|a'|\le N-1,\\
    &\int_0^t<t'>^{-1.99}E_N(t')E^2_{|a'|+1}(t')dt',\quad|a'|\le N-3.
\end{aligned}
\right.
\end{array}
\end{equation*}
}

{\bf Proof.} It is easy to verify that
\begin{equation}\label{g:ineq'}
    \|\Gamma^ag(t,r)\|^2_{L^2}
    \ls\|\Gamma^ag(0,r)\|^2_{L^2}+\left|\int_0^t\int\Gamma^ag\p_t\Gamma^agrdrdt'\right|
    \ls\|\Gamma^ag(0,r)\|^2_{L^2}+J_1+J_2+J_3,
\end{equation}
where
\begin{eqnarray*}
    &&J_1=:\left|\int_0^t\int\Gamma^ag\Gamma^af(\p_rg+\frac1rg)rdrdt'\right|,\\
    &&J_2=:\left|\int_0^t\int\Gamma^agf(\p_r+\frac1r)\Gamma^agrdrdt'\right|,\\
    &&J_3=:\sum_{b+c=a,b<a,c<a}\int_0^t\int\left|\Gamma^ag\frac{\Gamma^bf}{r}(r\p_r+1)\Gamma^cg\right|rdrdt'.
\end{eqnarray*}

By the same analysis to the first term $\Gamma^af\Gamma^a(g^2)$ as in $I_{04}$ of \eqref{I0:decomposition},
we achieve
\begin{equation}\label{J1:ineq}
    J_1\ls E_{\max\{3,[\frac{|a|}{2}]+2\}}(t)E^2_{|a|}(t)+\cQ_a(t).
\end{equation}

For $J_2$, we have
\begin{equation}\label{J2:ineq}
    J_2\ls\int_0^t\int|\Gamma^ag|^2\big[|\chi_0\p_rf|+\frac1r|\chi_0f|\big]rdrdt'
    \ls\int_0^t<t'>^{-1.99}E_3(t')E^2_{|a|}(t')dt'.
\end{equation}

For $J_3$, by $|(r\p_r+1)\Gamma^cg|\ls|(S+1)\Gamma^cg|+t|\p_t\Gamma^cg|$,
we then obtain
\begin{equation}\label{J3:ineq}
\begin{array}{l}
\displaystyle J_3\ls\sum_{b+c=a,|b|\le N-3}
\int_0^t<t'>^{-1.99}E^2_{|a|}(t')E_{|b|+3}(t')dt'\\
\displaystyle \qquad+\left\{\begin{aligned}
    &\int_0^t<t'>^{-1}E_{N-2}(t')E^2_{|a|}(t')dt',\quad|a|\le N,\\
    &\int_0^t<t'>^{-1.99}E_{|a|+2}(t')E^2_{|a|}(t')dt',\quad|a|\le N-2,
\end{aligned}
\right.\\
\displaystyle \ls\cQ_a(t).
\end{array}
\end{equation}

Collecting \eqref{g:ineq'}--\eqref{J3:ineq} yields \eqref{g:ineq}.

Next, we treat the case of $|a'|\le N-1$ in \eqref{curl:L2norm}. Similarly to the proof of \eqref{g:ineq}, we have
\begin{equation}\label{drg:ineq'}
\begin{array}{l}
\displaystyle \|\Gamma^{a'}w(t,r)\|^2_{L^2}\ls\|\Gamma^{a'}w(0,r)\|^2_{L^2}
+\left|\int_0^t\int\Gamma^{a'}w\p_t\Gamma^{a'}wrdrdt'\right|\\
\displaystyle \qquad\ls\|\Gamma^{a'}w(0,r)\|^2_{L^2}+J'_1+J'_2,
\end{array}
\end{equation}
where
\begin{eqnarray*}
    &&J'_1=:\left|\int_0^t\int\Gamma^{a'}wf\p_r\Gamma^{a'}wrdrdt'\right|,\\
    &&J'_2=:\sum_{c+d=a',d<a'}\int_0^t\int\left|\Gamma^{a'}w\Gamma^cf\p_r\Gamma^dw\right|rdrdt'.
\end{eqnarray*}

It is easy to check that $\|\Gamma^{a'}w\|_{L^2}=\|\Gamma^{a'}[(1+v)(\p_rg+\frac1rg)]\|_{L^2}\ls E_{|a'|+1}(t)$.
Thus, by integrating with respect to $r$, one has that for $J'_1$,
\begin{equation}\label{J1':ineq}
    J'_1\ls\int_0^t\int|\Gamma^{a'}w|^2\big[|\chi_0\p_rf|+\frac1r|\chi_0f|\big]rdrdt'
    \ls\int_0^t<t'>^{-1.99}E_3(t')E^2_{|a'|+1}(t')dt'.
\end{equation}

With the help of \eqref{dr:regular}, we easily get that
\begin{eqnarray*}
    |r\p_r\Gamma^dw|\ls|S\Gamma^dw|+t|\p_t\Gamma^dw|\ls|S\Gamma^dw|
    +t\sum_{b'+c'=d}|\frac1r\chi_0\Gamma^{b'}f||r\p_r\Gamma^{c'}w|\\
    \ls|S\Gamma^dw|+M\ve_0\sum_{c'\le d}|r\p_r\Gamma^{c'}w|,
\end{eqnarray*}
which implies
\begin{equation}\label{rdr:curl}
    |r\p_r\Gamma^dw|\ls|S\Gamma^dw|.
\end{equation}
Then, for the higher-order derivatives with $N-3\le|a'|\le N-1$, we achieve
\begin{equation}\label{J2':higher}
\begin{array}{l}
\displaystyle J'_2\ls\sum_{|c|\le N-3}\int_0^t<t'>^{-1.99}E_{|c|+3}(t')E^2_{|a'|+1}(t')dt'
+\sum_{|d|\le N-5}\int_0^t<t'>^{-1}\|S\Gamma^dw\|_{L^\infty}E^2_{|a'|+1}(t')dt'\\
\displaystyle \ls\int_0^t<t'>^{-1.99}E_N(t')E^2_{|a'|+1}(t')dt'+\int_0^t<t'>^{-1}E_{N-2}(t')E^2_{|a'|+1}(t')dt'.
\end{array}
\end{equation}
While, for the lower-order derivatives with $|a'|\le N-3$, we have
\begin{equation}\label{J2':lower}
    J'_2\ls\sum_{|c|\le|a'|\le N-3}\int_0^t<t'>^{-1.99}E_{|c|+3}(t')E^2_{|a'|+1}(t')dt'
    \ls\int_0^t<t'>^{-1.99}E_N(t')E^2_{|a'|+1}(t')dt'.
\end{equation}

Substituting \eqref{J1':ineq}--\eqref{J2':lower} into \eqref{drg:ineq'} yields \eqref{curl:L2norm}.

Next, we treat the case of $|b|\le N-3$ in \eqref{curl:L3norm}.
Note that
\begin{equation*}
    \|\Gamma^bw(t,r)\|^3_{L^3}\ls\|\Gamma^bw(0,r)\|^3_{L^3}
    +\left|\int_0^t\int|\Gamma^bw|\Gamma^bw\p_t\Gamma^bwrdrdt'\right|.
\end{equation*}
According to the equation of the specific vorticity \eqref{curl:eqn}, we have
\begin{eqnarray*}
    \left|\int_0^t\int|\Gamma^bw|\Gamma^bw\p_t\Gamma^bwrdrdt'\right|
    \ls\int_0^t\int|\Gamma^bw|^3|\p_r(rf)|drdt'\\ +\sum_{c+d=b,c<b}\int_0^t\int|\Gamma^bw|^2|r\p_r\Gamma^cw|\frac{|\chi_0\Gamma^df|}{r}rdrdt'.
\end{eqnarray*}
Applying \eqref{rdr:curl} again, we obtain
\begin{equation*}
    \|\Gamma^bw(t,r)\|^3_{L^3}\ls\|\Gamma^bw(0,r)\|^3_{L^3}
    +\int_0^t<t'>^{-1.99}E_{|b|+3}(t')W^3_{|b|}(t')dt'.
\end{equation*}
Then, \eqref{curl:L3norm} is proved.

Similarly, for $|b|\le N-3$ in \eqref{curl:W13norm}, we arrive at
\begin{equation*}
    \|\p_r\Gamma^bw(t,r)\|^3_{L^3}\ls\|\p_r\Gamma^bw(0,r)\|^3_{L^3}
    +\left|\int_0^t\int|\p_r\Gamma^bw|\p_r\Gamma^bw\p_r\p_t\Gamma^bwrdrdt'\right|.
\end{equation*}
By using \eqref{curl:eqn} again, we get
\begin{equation*}
\begin{array}{l}
\displaystyle \left|\int_0^t\int|\p_r\Gamma^bw|\p_r\Gamma^bw\p_r\p_t\Gamma^bwrdrdt'\right| \ls\int_0^t\int|\p_r\Gamma^bw|^3|\p_r(rf)|drdt'\\
\displaystyle \qquad+\sum_{c+d=b,c<b}\int_0^t\int|\p_r\Gamma^bw|^2|r\p_r\p_r\Gamma^cw|\frac{|\chi_0\Gamma^df|}{r}rdrdt'\\
\displaystyle \qquad+\sum_{c+d=b}\int_0^t\int|\p_r\Gamma^bw|^2|\p_r\Gamma^cw||\chi_0\p_r\Gamma^df|rdrdt'.\\
\end{array}
\end{equation*}
Analogously, thanks to \eqref{dr:regular}, we deduce that
\begin{eqnarray*}
    &&|r\p_r\p_r\Gamma^cw|\ls|\p_rS\Gamma^cw|+|\p_r\Gamma^cw|+t'|\p_t\p_r\Gamma^cw|\\
    &&\ls|\p_rS\Gamma^cw|+|\p_r\Gamma^cw|+t'\sum_{d'+e=c}\big[|\chi_0\p_r\Gamma^{d'}f||\p_r\Gamma^ew|
    +|\frac1r\chi_0\Gamma^{d'}f||r\p_r\p_r\Gamma^ew|\big]\\
    &&\ls|\p_rS\Gamma^cw|+\sum_{e\le c}|\p_r\Gamma^ew|
    +M\ve_0\sum_{e\le c}|r\p_r\p_r\Gamma^ew|,
\end{eqnarray*}
which implies
\begin{equation*}
    \|\p_r\Gamma^bw(t,r)\|^3_{L^3}\ls\|\p_r\Gamma^bw(0,r)\|^3_{L^3}
    +\int_0^t<t'>^{-1.99}E_{|b|+3}(t')W^3_{|b|}(t')dt'.
\end{equation*}
Thus, \eqref{curl:W13norm} is achieved. \hfill $\square$

\section{Proof of Theorem 1.1}\label{sect8}

Based on the results in Section 3--Section 7, we now start to prove Theorem 1.1.
\vskip 0.1 true cm

{\bf Proof of Theorem~1.1.} For $|a|\le N-1$, it holds that
\begin{equation*}
    \|(\p_r+\frac1r)\Gamma^ag(t,r)\|_{L^2}=\|\tilde\Gamma^a(\p_r+\frac1r)g(t,r)\|_{L^2}
    =\|\tilde\Gamma^a\Big(\frac{w}{1+\vv+G}\Big)\|_{L^2}.
\end{equation*}
Then, combining this with \eqref{energyineq}, \eqref{nearLC:ineq}, \eqref{awayLC:ineq},  Lemma~7.2 and
assumption \eqref{assumption}, we achieve
\begin{eqnarray*}
    E^2_N(t)\le CE^2_N(0)+CM\ve\int_0^t<t'>^{-1}E^2_N(t')dt',\\
    E^2_{N-2}(t)\le CE^2_{N-2}(0)+CM\ve\int_0^t<t'>^{-\frac54}E^2_{N-2}(t')dt',\\
    W^3_{N-3}(t)\le CW^3_{N-3}(0)+CM\ve\int_0^t<t'>^{-\frac54}W^3_{N-3}(t')dt'.
\end{eqnarray*}

It is not hard to check that $E^2_N(0)+W^3_{N-3}(0)\le C\ve$.
Applying the Gronwall's inequality to the above inequalities and then using \eqref{WT:energy1} and \eqref{WT:energy0},
we deduce that for sufficiently large but fixed positive constants $M,M'$, and for small $\ve_0>0$ enough,
\begin{eqnarray*}
    E_N(t)+\cX_N(t)+\cY_N(t)\le C\ve(1+t)^{CM\ve}\le\frac12M\ve(1+t)^{M'\ve},\\
    E_{N-2}(t)+\cX_{N-2}(t)+\cY_{N-2}(t)+W_{N-3}(t)\le C\ve e^{CM\ve}\le\frac12M\ve.
\end{eqnarray*}
Therefore, Theorem~1.1 is achieved by the local existence of smooth solution
to \eqref{EulerC1form} with \eqref{ChaplyginGas} and the continuity argument. \hfill $\square$


\begin{thebibliography}{99}

\bibitem{Alinhac92} S.~Alinhac, {\it Une solution approch\'ee en grand temps des \'equations d'Euler compressibles
axisym\'etriques en dimension deux,}
 Comm. Partial Differential Equations \textbf{17} (1992), no. 3-4, 447--490.

\bibitem{Alinhac93} S. Alinhac, {\it Temps de vie des solutions r\'eguli\'eres des \'equations d'Euler compressibles axisym\'etriques en dimension deux,}
  Invent. Math. \textbf{111} (1993), 627--670.

\bibitem{Alinhac99} S. Alinhac, {\it Blowup of small data solutions for a class of quasilinear
wave equations in two space dimensions. II},
Acta Math. \textbf{182} (1999), no. 1, 1--23.

\bibitem{Alinhac01} S. Alinhac, {\it The null condition for quasilinear wave equations in two space dimensions I,}
Invent. Math. \textbf{145} (2001), no. 3, 597--618.


\bibitem{CF:book} R. Courant, K. O. Friedrichs, {\it Supersonic flow
and shock waves, Interscience Publishers Inc., New York}, 1948.



\bibitem{Christodoulou86} D. Christodoulou, \textit{Global solutions of nonlinear
hyperbolic equations for small initial data}, Comm. Pure Appl.
Math. \textbf{39} (1986), no. 2, 267--282.

\bibitem{Christodoulou07} D. Christodoulou, \textit{The formation of shocks in 3-dimensional fluids,} EMS Monogr. Math., Eur. Math. Soc., Z\"urich, 2007.

\bibitem{CM} D. Christodoulou, Miao Shuang, \textit{Compressible flow and Euler's equations}, Surveys of Modern Mathematics, 9. International Press, Somerville, MA; Higher Education Press, Beijing, 2014.



\bibitem{DWY15} Bingbing Ding, Ingo Witt, Huicheng Yin, {\it The global smooth symmetric solution to 2-D full
compressible Euler system of Chaplygin gases,} J. Differential Equations \textbf{258} (2015), no. 2, 445--482.


\bibitem{DLY15} Bingbing Ding, Yingbo Liu, Huicheng Yin,  \textit{The small data solutions of general 3D quasilinear wave equations. I,} SIAM J. Math. Anal. \textbf{47} (2015), no. 6, 4192--4228.

\bibitem{DWY16} Bingbing Ding, Ingo Witt, Huicheng Yin, \textit{The small data solutions of general 3-D quasilinear wave equations. II,} J. Differential Equations \textbf{261} (2016), no. 2, 1429--1471.

\bibitem{Godin05} P. Godin, {\it The lifespan of a class of smooth spherically symmetric solutions of the compressible Euler equations with variable entropy in three space dimensions,} Arch. Ration. Mech. Anal. \textbf{177} (2005), no. 3, 479--511.


\bibitem{Godin07} P. Godin, {\it Global existence of a class of smooth 3D spherically symmetric flows of Chaplygin
gases with variable entropy,} J. Math. Pures Appl. \textbf{87} (2007), 91--117.



\bibitem{HKSW} G. Holzegel, S. Klainerman, J. Speck, W.W.-Y. Wong, {\it Small-data shock formation in solutions to 3d quasilinear wave equations: An overview,} Journal of Hyperbolic Differential Equations \textbf{13} (2016), no. 01, 1--105.




\bibitem{Hormander97book} L. H\"ormander, \textit{Lectures on nonlinear hyperbolic differential equations.} Math\'ematiques \& Applications (Berlin) [Mathematics \& Applications], 26. Springer-Verlag, Berlin, 1997. viii+289 pp.

\bibitem{John} F. John, \textit{Nonlinear wave equations, formation of singularities.}
    Seventh Annual Pitcher Lectures delivered at Lehigh University, Bethlehem, Pennsylvania, April 1989. University Lecture Series, 2. American Mathematical Society, Providence, RI, 1990.

\bibitem{Klainerman} S. Klainerman, {\it The null condition and global existence to nonlinear wave equations,} in: Nonlinear Systems of Partial Differential Equations in Applied Mathematics, Part \textbf{1}, Santa Fe, NM, 1984, in: Lect. Appl. Math., vol. \textbf{23}, Amer. Math. Soc., Providence, RI, 1986, pp. 293--326.

\bibitem{Lei16} Zhen Lei, {\it Global well-posedness of incompressible elastodynamics in two dimensions,} Comm. Pure Appl. Math. \textbf{69} (2016), no. 11, 2072--2106.

\bibitem{LeiWei17} Zhen Lei, Changhua Wei, {\it Global radial solutions to 3D relativistic Euler equations for non-isentropic Chaplygin gases,} Math. Ann. \textbf{367} (2017), no. 3-4, 1363--1401.

\bibitem{LiYin18} Jun Li, Huicheng Yin, {\it Global smooth solutions of 3-D null-form wave equations in exterior domains with Neumann boundary conditions,} J. Differential Equations \textbf{264} (2018), no. 2, 5577--5628.

\bibitem{Hans0} H. Lindblad, {\it On the lifespan of solutions of nonlinear wave equations with small initial data,} Comm. Pure Appl.
Math. \textbf{43} (1990), no. 4, 445--472.


\bibitem{Lindblad08AJM} H. Lindblad, {\it Global solutions of quasilinear wave equations,} Amer. J. Math. \textbf{130} (2008), no. 1, 115--157.

\bibitem{LNS13} H. Lindblad, M. Nakamura, C.D. Sogge, {\it Remarks on global solutions for nonlinear wave equations under the standard null conditions,} J. Differential Equations \textbf{254} (2013), no. 3, 1396--1436.

\bibitem{LukSpeck16} J. Luk, J. Speck, {\it Shock formation in solutions to the 2D compressible Euler equations in the presence of non-zero vorticity,} arXiv:1610.00737 (2016)

%\bibitem{Majda:book} A. Majda, {\it Compressible fluid flow and systems of conservation laws in several space variables,} Applied Mathematical Sciences, 53. Springer-Verlag, New York, 1984.


\bibitem{MS04} A. Morando, P. Secchi, {\it On 3D slightly compressible Euler equations,}
 Port. Math. (N.S.) \textbf{61} (2004), no. 3, 301--316.


\bibitem{Rammaha89} M.~Rammaha, \textit{Formation of singularities in compressible fluids in two-space dimensions,}
  Proc. Am. Math. Soc. \textbf{107} (1989), 705--714.


\bibitem{Secchi02} P. Secchi, {\it On slightly compressible ideal flow in the half-plane,}
Arch. Ration. Mech. Anal. \textbf{161} (2002), no. 3, 231--255.


\bibitem{Sideris85} T. Sideris, \textit{Formation of singularities in three-dimensional compressible fluids,} Comm. Math. Phys. \textbf{101} (1985), 475--485.

\bibitem{Sideris97} T. Sideris, {\it Delayed singularity formation in 2D compressible flow,} Amer. J. Math. \textbf{119} (1997), 371--422.

%\bibitem{Sideris00} T. Sideris, {\it Nonresonance and global existence of prestressed nonlinear elastic waves,} Ann. of Math. (2) \textbf{151} (2000), no. 2, 849--874.

\bibitem{Speck17} J. Speck, {\it A summary of some new results on the formation of shocks
in the presence of vorticity.} Nonlinear analysis in geometry and applied mathematics, 133-157,
Harv. Univ. Cent. Math. Sci. Appl. Ser. Math., 1, Int. Press, Somerville, MA, 2017.


\bibitem{Speck16book} J. Speck, {\it Shock formation in small-data solutions to 3D quasilinear wave equations.}
Mathematical Surveys and Monographs, 214. American Mathematical Society, Providence, RI, 2016. xxiii+515 pp.

\bibitem{Yin} Huicheng Yin, \textit{Formation and construction of a shock wave for 3-D compressible
Euler equations with the spherical initial data,}
  Nagoya Math. J., \textbf{175} (2004), 125--164.







\end{thebibliography}
\end{document}